\documentclass{article}
\usepackage{tabularx}
\usepackage{graphicx}
\usepackage{adjustbox}
\usepackage[utf8]{inputenc}
\usepackage[utf8]{inputenc}
\usepackage{multirow}
\usepackage{longtable,verbatim}
\usepackage{enumitem,  rotfloat, supertabular}
\usepackage{booktabs}
\usepackage{mathtools}
\usepackage{amsmath,hhline,dcolumn,layout,booktabs, pdflscape, amsthm}
\usepackage{amsfonts}
\usepackage{todonotes}
\usepackage{algpseudocode}
\usepackage{url}
\allowdisplaybreaks

\floatstyle{ruled}
\floatname{algorithm}{Algorithm}

\usepackage[a4paper,top=1in,bottom=1in,left=1in,right=1in]{geometry}
\usepackage[utf8]{inputenc}
\usepackage[english]{babel}
\usepackage{natbib}
\usepackage{pgfplots}
\usepackage{tikz}
\usetikzlibrary{patterns}
\usepackage{url}
\setcitestyle{authoryear,open={(},close={)}}

\definecolor{ggreen}{HTML}{9BBB59}
\definecolor{rred}{HTML}{C0504D}
\usepackage{graphicx}
\usepackage{subcaption}

\newcases{nullcases}
    {\quad}
    {\hfil{##}} {{##}\hfil}
    {.} {.}

\makeatother

\usepackage{hyperref}
    \hypersetup{
        colorlinks   = true,
        citecolor    = blue,
        linkcolor=blue  
    }

\title{Capacity Planning for Effective Cohorting of Dialysis Patients during the Coronavirus Pandemic: A Case Study}
\author{Cem D.C. Bozkir $^{a}$ $\bullet$ Cagri Ozmemis $^{a}$ $\bullet$ Ali Kaan Kurbanzade $^{a}$ \\  Burcu Balcik $^{a,\ast}$ $\bullet$ Evrim D. Gunes $^{b}$ $\bullet$ Serhan Tuglular $^{c}$\\
        \small $^{a}$ Industrial Engineering Department, Ozyegin University, Istanbul, Turkey \\
        \small $^{b}$Business Administration, College of Administrative Sciences
and Economics, \\ \small Koc University, Sariyer, Istanbul, Turkey \\
        \small $^{c}$Medical Faculty, Department of Internal Medicine, Marmara University, Istanbul, Turkey \\\\
        \small $^{*}$Corresponding author: \tt{burcu.balcik@ozyegin.edu.tr} \\
}

\date{}

\begin{document}
\maketitle

\normalsize
\begin{abstract}
Chronic dialysis patients have been among the most vulnerable groups of the society during the coronavirus (COVID-19) pandemic as they need regular treatments in a hospital environment, facing infection risk. Moreover, the demand for dialysis resources has significantly increased since many COVID-19 patients need acute dialysis due to kidney failure. In this study, we address capacity planning decisions of a hemodialysis clinic located within a major hospital in Istanbul, designated to serve both infected and uninfected patients during the pandemic with limited resources (i.e., dialysis machines). The hemodialysis clinic applies a three-unit cohorting strategy to treat four types of patients in separate units and at different times to mitigate infection spread risk among patients. Accordingly, at the beginning of each week, the clinic needs to determine the number of available dialysis machines to allocate to each unit that serves different patient cohorts. Given the uncertainties in the number of different types of patients that will need dialysis, it is a challenge to allocate the scarce dialysis resources effectively by evaluating which capacity configuration would minimize the overlapping treatment sessions of infected and uninfected patients over a week. We represent the uncertainties in the number of patients by a set of scenarios and present a two-stage stochastic programming approach to support capacity allocation decisions of the clinic. We present a case study based on the real-world patient data obtained from the clinic to illustrate the effectiveness of the proposed modeling approach and compare the performance of different cohorting strategies. \\

\end{abstract}

\newpage
\section{Introduction}

The coronavirus disease (COVID-19) has caused a global pandemic,  as declared by the World Health Organization in March 2020 \citep{sohrabi2020world}. It has created a significant burden on the capacity of healthcare systems all around the world due to a surge in demand for outpatient clinics, inpatient beds, and intensive care units. While many people could postpone their nonurgent hospital visits to protect from infection risk \citep{nyt}, some patients with care needs due to chronic health conditions (such as dialysis, chemotherapy, physical therapy) do not have the luxury to avoid or postpone visits to health care facilities since they are obliged to obtain medical care regularly, making them  more vulnerable. It has also been a challenge for the hospitals to serve such chronic patients without exposing them to infection risk with ever-tighter resources. This paper addresses the capacity planning decisions of a hospital’s hemodialysis clinic that provides care to infected and uninfected patients during the COVID-19 pandemic by implementing cohorting strategies to mitigate the risk of infection spread among patients.


Dialysis is a treatment to remove waste products and excess fluid from the blood with the help of a machine when the kidneys fail to work \citep{nhs}. Most chronic hemodialysis patients must receive dialysis treatment three times a week. During the treatment, patients stay connected to a dialysis machine for about four hours, which filters and purifies the blood. During their frequent and long hospital visits to get treatments, chronic hemodialysis patients can become highly susceptible to COVID-19 transmission risk. Indeed, their infection risk is found as nearly two-fold greater than those receiving dialysis at home \citep{hsu2020covid}. Moreover, Chronic Kidney Disease (CKD) is the most prevalent risk factor for severe COVID-19 cases \citep{era2021chronic}, 
and mortality rate of patients with CKD under COVID-19 virus infection is significantly higher compared to those without CKD \citep{roper2020delivering,ozturk2020mortality}. Therefore, it is imperative for health care providers to take precautions to ensure that the uninfected chronic dialysis patients are not in contact with the COVID-19 infected patients during their clinic visits \citep{basile2020recommendations}. Given that the global prevalence rate of chronic kidney disease (all stages) is 9.1\% \citep{bikbov2020}, effective planning of dialysis treatments during the pandemics by accounting for infection spread risk can protect the lives of millions of people in vulnerable conditions. However, even without a pandemic disturbance, planning and scheduling of treatments in a single dialysis unit is already a difficult problem (e.g., \cite{liu2019patient,holland1994scheduling,fleming2019analytical}). When resources become tighter, such as the aftermath of an earthquake \citep{sever2009} or during a pandemic \citep{corbett2020epidemiology}, the increased demands and uncertainties make it ever more challenging to plan for dialysis services.

In this study, motivated by a real case, we focus on dialysis services provided in a hospital during the pandemics where multiple patient types must be served. Specifically, besides the chronic hemodialysis patients, who need regular care, some patients need dialysis due to acute kidney failure for various medical conditions. These acute patients, who often have multiple health complications, are already admitted to the hospital, and transferring them to another hospital is usually impossible. During the pandemic, the hospital is additionally responsible for treating the suspected and infected COVID-19 patients with dialysis needs that cannot be treated in unequipped small dialysis centers. We address the capacity planning decisions of the hospital to treat different patient types by mitigating the transmission risk between infected and uninfected patients.

Managing dialysis resources effectively during the pandemic is challenging for the hospital due to several reasons. First, the effective dialysis capacity can be reduced due to the additional time required for cleaning the rooms in-between dialysis sessions as well as arrangements required to cohort patients. Specifically, the uninfected (clean), infected (confirmed), and suspected (shows symptoms or had contact with a confirmed case) patients must be treated at different designated units; and if possible, during different (nonoverlapping) periods. Secondly, the demand for dialysis during a pandemic increases because acute dialysis sessions may be needed for some infected patients due to the adverse effects of COVID-19 on kidneys \citep{sperati2020}; indeed, 35.2\% of COVID-19 patients needed dialysis \citep{klein2020, smith2020}. Thirdly, while the number of chronic patients on a treatment schedule is known, the number of infected and clean acute patients that need dialysis treatments each day can be highly uncertain. In this resource-constrained environment, employing an operations research based approach can be valuable to make effective cohorting, capacity planning and treatment scheduling decisions by using limited resources efficiently.

Treating the infected, suspected, and uninfected patients at different times and in separate areas are common approaches followed in healthcare facilities to mitigate the risk of infection transmission among patients \citep{kliger2020mitigating}. Different cohorting strategies, which involve defining patient types and assigning them to separate areas, have been implemented by health care providers around the world during the COVID-19 pandemic (e.g., \citet{collison2020, whiteside2020}), and there is no well-established standard yet. \citet{park2020korean} and \citet{meijers2020safeguarding} summarize the classification of different hemodialysis patients that are either infected by COVID-19 or under the suspicion of infection and outline the operationalization of the dialysis clinics in Korea, Belgium, and Italy. Despite its importance, no study has yet proposed analytical tools to support the effective implementation of cohorting strategies in a dialysis clinic during the pandemic.

In this paper, we focus on the operations of a dialysis clinic in a large public teaching and research hospital located in Istanbul during the COVID-19 pandemic. This hospital is one of the designated ``pandemic hospitals" in this large densely-populated city. The dialysis clinic is responsible for both serving the uninfected chronic dialysis patients and acute inpatients, as well as the infected and suspected dialysis patients, 
often through referrals from other clinics/hospitals that cannot treat COVID-19 patients. The clinic has been implementing a three-tier (unit) cohorting strategy since the beginning of the pandemic to mitigate infection spread risk among patients. To facilitate this, the clinic is divided into units separated by drywalls, where each unit has several machines. Accordingly, uninfected, suspected and infected patients are treated at different units.  Depending on the estimated demands for the upcoming week, the hospital assigns dialysis machines to the units in the dialysis clinic thereby setting their capacity to serve different patient cohorts during the week. Moreover, it is desirable to schedule the dialysis treatments of different patients at different times during the day to minimize physical interaction and avoid infection spread. An alternative cohorting strategy could separate the clinic into two units \citep{whiteside2020}, one for the uninfected patients, and the other for suspected and infected patients, where these patients would be treated sequentially in the unit.

Our collaboration with the dialysis clinic has two main objectives: (1) To develop mathematical models for making capacity planning decisions to implement a given cohorting strategy effectively, (2) To evaluate the performance of alternative cohorting strategies with two and three tiers (units). To achieve these objectives, we develop mathematical models to solve the capacity planning 
associated with each cohorting strategy. Specifically, for a given (i.e., two- or three-unit) cohorting strategy, we represent the uncertainty in the number of dialysis patients in each cohort that will need dialysis for the upcoming week by scenarios and develop two-stage stochastic programming models to decide on the number of dialysis machines that must be allocated to each unit at the beginning of a week (first stage), based on the second stage decisions which determine a daily dialysis treatment schedule for each unit by considering the number of patients in a scenario. The objective is to treat all patients in the clinic while minimizing the expected number of patients from different cohorts having dialysis at the overlapping sessions. 
We test the proposed models on real-world data obtained from our collaborating hospital regarding demands for different patient types and cohorting plans 
over a two-month period during the pandemic. We present results that illustrate the benefits of using the proposed models and compare the performance of different cohorting strategies.

The remainder of this paper is structured as follows. In \S \ref{s:litreview}, we review the related literature. In \S \ref{s:problem}, we describe the system, define the problem and present our models. In \S \ref{s:case}, we present a case study and discuss the results. Finally, we present our conclusions in \S \ref{s:conclusion}.

\section{Literature Review} \label{s:litreview}

We review the literature related to health care capacity planning problems in resource-constrained settings, and planning of hemodialysis services during pandemics.

\subsection{Health care capacity planning in resource-limited settings}

There exists abundant literature that utilizes operations research methods to address capacity planning and resource allocation problems in health care settings. Existing studies address problems in different settings ranging from planning the operations of a single hospital or a unit (e.g., intensive care units, operating rooms) to a network of health providers (see reviews by \citet{guerriero2011, rais2011, hulshof2012, ahmadi2017, bai2018}). 

Catastrophic events (such as pandemics, natural and man-made disasters) can create unexpected demand surge and extreme resource scarcity in healthcare systems. Several studies address planning problems faced in such settings and propose various strategies to manage the demand surge by using available resources efficiently. For example, \citet{repoussis2016} focus on regional emergency planning and present a mixed integer programming (MIP) model for the combined ambulance dispatching, patient-to-hospital assignment, and treatment ordering problem. The proposed model considers the triage levels of casualties at each site and minimizes the maximum completion time of treatments. \citet{caunhye2018stochastic} present a scenario-based stochastic programming model to determine the location of the alternative care facilities to be established following a disaster, and the number of patients to be directed to these centers by integrating triage and self movement of casualties. \citet{ouyang2020allocation} focus on bed allocation decisions in an intensive care unit (ICU) of a hospital during times of high demand and present a stochastic model to determine whether a new patient will be admitted to the ICU and which patients will be discharged from the ICU to minimize the long-run average mortality rate. \citet{mills2020surge} evaluate mitigation and response strategies to create surge capacity for emergency patients in a hospital by taking disposition actions for different types of patients with varying severity. The authors present an optimization model and perform a simulation study to evaluate the performance of strategies based on the characteristics of individual hospitals.

Several studies specifically focus on managing health care resources at the time of pandemics. \citet{arora2010} focus on satisfying the large amounts of need for antiviral medicines during an influenza pandemic. The authors formulate a two-stage model to first decide on the locations and amounts of supplies to be prepositioned in the stockpile to prevent the disease spread, and then make redistribution decisions with respect to transshipment and stockpiling costs in the treatment phase while considering the expected number of people affected. \citet{sun2014} address allocation of inpatients among multiple hospitals over several months during an influenza outbreak by considering different patient types with respect to their needs for various medical equipment. The authors present a multi-objective model that minimizes the total distance traveled by patients and the maximum distance of a patient to a hospital. \citet{liu2015} develop an integer programming (IP) model to allocate medical resources to be transported from suppliers to hospitals through a set of distribution centers to control the spread of influenza pandemic. The authors base their demand estimations on a susceptible-exposed-infected-recovered (SEIR) model and aim to minimize the total transportation cost while satisfying the dynamically changing demand for the resources. \citet{long2018spatial} present a two-stage model that handles the allocation of treatment units used in Ebola response across different geographical regions according to the course of the epidemic, where the number of infected people is predicted by an SEIR model. The authors present alternative approaches to allocate health resources including a heuristic, a greedy policy, a linear program, and a dynamic programming algorithm.

Since the declaration of the COVID-19 pandemics in March 2020, there has been growing research attention on the challenges faced in allocating and managing valuable health care resources during the pandemics. \citet{klein2020} present a review of studies that present models and tools to estimate demand surge for hospital capacity and resources (such as hospital beds, ventilators) during the pandemics. \citet{mehrotra2020model} present a stochastic MIP model to allocate and share a critical resource during the COVID-19 pandemic, the ventilators, among risk-averse states by considering demand uncertainties. The authors apply their model to a case study by considering FEMA as a coordinator that controls the allocation of ventilators among the states in the USA. \citet{parker2020optimal} study a demand and resource redistribution problem among multiple health care facilities. The authors develop deterministic and robust optimization models to determine the optimal demand and resource transfers among facilities to minimize the required surge capacity and resource shortage during a period of heightened demand. 



Although cohorting patients is a widely implemented method to prevent infection spread risk during pandemics, a limited number of studies present analytical models to support the implementation of cohorting strategies in health care settings. \citet{pinker2013} address capacity configuration decisions in a hospital with the existence of patients needing isolation due to infectious diseases. The authors present a stochastic optimization model with a revenue maximization objective to decide on the beds to be reserved for isolation and non-isolation patients, where the two groups have different requirements for spacing and admission criteria. \citet{chia2015} analyze the effectiveness of different cohorting strategies in a pediatric hospital for allocating beds among different patient types. \citet{wang2020} focus on allocating beds among different departments in a hospital in a setting with stochastic patient arrival and service times and propose an approach that combines a MIP and a simulation model to minimize the expected cost of rejection and waiting. \citet{melman2020balancing} develop a discrete event simulation model to evaluate alternative strategies for allocating limited resources effectively in a hospital to prioritize the surgeries of COVID-19 and uninfected patients. The authors apply their model in a hospital and show that even though elective surgeries were canceled to use major hospital capacity for COVID-19 treatment, it is better to prioritize such surgeries until a threshold to minimize total deaths. 

Different from the existing studies, we focus on weekly capacity planning decisions in a dialysis clinic to treat multiple patient types in a pandemic setting. To the best of our knowledge, no study has yet focused on modeling and evaluating alternative cohorting strategies where both rooms and treatment times of different patient cohorts must be separated to mitigate infection spread, and the capacity configuration decisions are made to minimize expected overlapping sessions.

\subsection{Hemodialysis planning during pandemics}

The literature that addresses hemodialysis planning during pandemics is limited to practical medical guidelines to manage this process effectively \citep{naicker2020novel}. For instance, \citet{ikizler2020minimizing}, \citet{kliger2020mitigating} and \citet{rubin2020finding} highlight the importance of cohorting, screening and disinfection in order to prevent infection in dialysis clinics. \citet{CMS2020} suggests that uninfected dialysis patients should wait for their sessions outside the hospital and there should be at least two meters among patients in dialysis units. \citet{roper2020delivering} point out high infection risks in the waiting rooms of the clinics. \citet{alberici2020management} describe the layout changes made in a hospital to provide treatments to uninfected and infected patients with different kidney diseases. \citet{kliger2020mitigating}, \citet{roper2020delivering} and \citet{corbett2020epidemiology} recommend treating suspected or confirmed COVID-19 cases in isolated rooms and uninfected dialysis patients in regular rooms.  \citet{rubin2020finding} cohort the dialysis patients into asymptomatic, suspected and confirmed COVID-19 cases. 

Although there exist a few operations research studies that address problems related to planning hemodialysis treatments, no study has yet addressed the particular challenges faced during a pandemic setting. \citet{liu2019patient} study a dialysis scheduling problem with different types of patients (conventional, hepatitis B, and hepatitis C) that must be treated in a dialysis center. The authors present an IP model to minimize night shifts and meet patient preferences. \citet{fleming2019analytical} address the problem of scheduling dialysis patients by considering each dialysis machine as a workstation. The proposed IP model minimizes the waiting time of patients arriving at the clinic for treatment and the scheduled finish time of treatments each day. \citet{yu2020} study appointment scheduling policies for patients that need a series of appointments such as chemotherapy or chronic dialysis patients. The authors present a Markov decision process model that takes into account revenues per service per patient, and the costs of staffing, overtime, overbooking, and delay.
 
\citet{currie2020simulation} express the need for scientific studies to prevent the infection rate of COVID-19 in hemodialysis units and highlight that this area is lacking in the literature. \citet{corbett2020epidemiology} discuss the challenges in planning the treatments of dialysis patients under a high level of uncertainty and stress the need for an analytical method. \citet{ikizler2020minimizing} highlight that the prolonged pandemic overextends the healthcare capacities and creates shortages in equipment, which necessitates allocating resources most efficiently by cohorting symptomatic and clean patients. Although several studies draw attention to the need for effective management of dialysis resources during a pandemic, no work in the literature provides analytical methods for supporting the decision makers to make cohorting and the capacity planning decisions in a dialysis unit to serve different types of patients. We contribute to the literature by introducing a new problem and presenting mathematical models to support decision making for configuring the clinic's capacity to effectively plan dialysis sessions during the pandemic. We focus on four types of patients and present models for alternative cohorting strategies that divide the clinic into two and three units. The proposed models can be extended to other settings that consider different cohorting strategies with different number of cohorts and patient types.  



\section{Problem Description \& Mathematical Modeling} \label{s:problem}
In this section, we first describe the system based on the operations of our collaborating hemodialysis clinic located in a major hospital in Turkey (\S \ref{s:system}). We then define the problem (\S \ref{s:prob}) and present our formulations developed for two alternative cohorting strategies (\S \ref{s:model}).

\subsection {System Description and Objectives of the Study} \label{s:system}

We focus on the operations of a hemodialysis clinic that treats different types of patients during the COVID-19 pandemic. The clinic is located within the Marmara University's Pendik Training and Research Hospital, which is a major public hospital located in Istanbul. In Turkey, all major public hospitals have been declared as ``pandemic hospitals'' after COVID-19 was officially announced to be a pandemic in March 2020. A pandemic hospital is a designated hospital for treating COVID-19 infected or suspected patients. Nevertheless, the pandemic hospitals have continued providing their regular care services as well, which. necessitated using resources efficiently more than ever. While all COVID-19 patients could be treated only at the pandemic hospitals during the initial months of the pandemic, other hospitals started accepting COVID-19 patients in the later stages of the pandemic.

In Turkey, chronic hemodialysis patients can get their treatment either in a public hospital or a private dialysis clinic. The dialysis treatments are covered by state-wide health insurance. However, some patients that are in special status (such as refugees) can receive treatments only at the dialysis clinics located within the public hospitals. All chronic dialysis patients receive treatments at their registered clinics three days a week, which can be either on the Monday-Wednesday-Friday (MWF) or the Tuesday-Thursday-Friday (TTF) regime. Therefore, each clinic makes a schedule to treat its chronic patients by following these regimes. All dialysis clinics in the country handle the transportation of the patients to and from the clinic to avoid delays in the schedule due to possible late arrivals of the patients. Different from the small private dialysis centers, major public hospitals also have inpatients (i.e., admitted to various wards, intensive care units) who may need temporary dialysis treatment during their stay at the hospital (e.g., after an operation). The number of such acute patients that need dialysis treatment varies and the hospital arranges the treatments of the acute patients along with the chronic patients each day. However, the COVID-19 pandemic required making significant changes in this routine in our case hospital, which has additionally become responsible for serving the COVID-19 patients with the same level of resources (dialysis machines, rooms). 

 The private dialysis centers are only specialized in dialysis service and they are not fully equipped for managing possible complications of dialysis patients with COVID-19. 
Therefore, chronic dialysis patients, who were normally treated in such dialysis centers, are referred to pandemic hospitals if they become infected with COVID-19. These patients are treated at the hospital until they are fully recovered from COVID-19. Additionally, some COVID-19 patients who do not have a priori chronic kidney disease may need dialysis due to COVID-19 related complications. These patients are also treated at pandemic hospitals. Our case hospital's dialysis clinic has to admit all confirmed and suspected COVID-19  patients for dialysis treatment. The number of chronic or acute COVID-19 patients that will arrive at the dialysis clinic each day is not known beforehand. Hence, it is important for the hospital to make effective capacity planning to accommodate the treatments of the COVID-19 patients along with its regular patients under uncertainty. In this study, we focus on the cohorting strategies applied in the hospital to treat infected, suspected and uninfected patients to minimize physical interaction among different groups.

\begin{figure}[!hbt]
\begin{center}
\caption{Dialysis Clinic Floor Plan} \medskip
\includegraphics[scale=0.35]{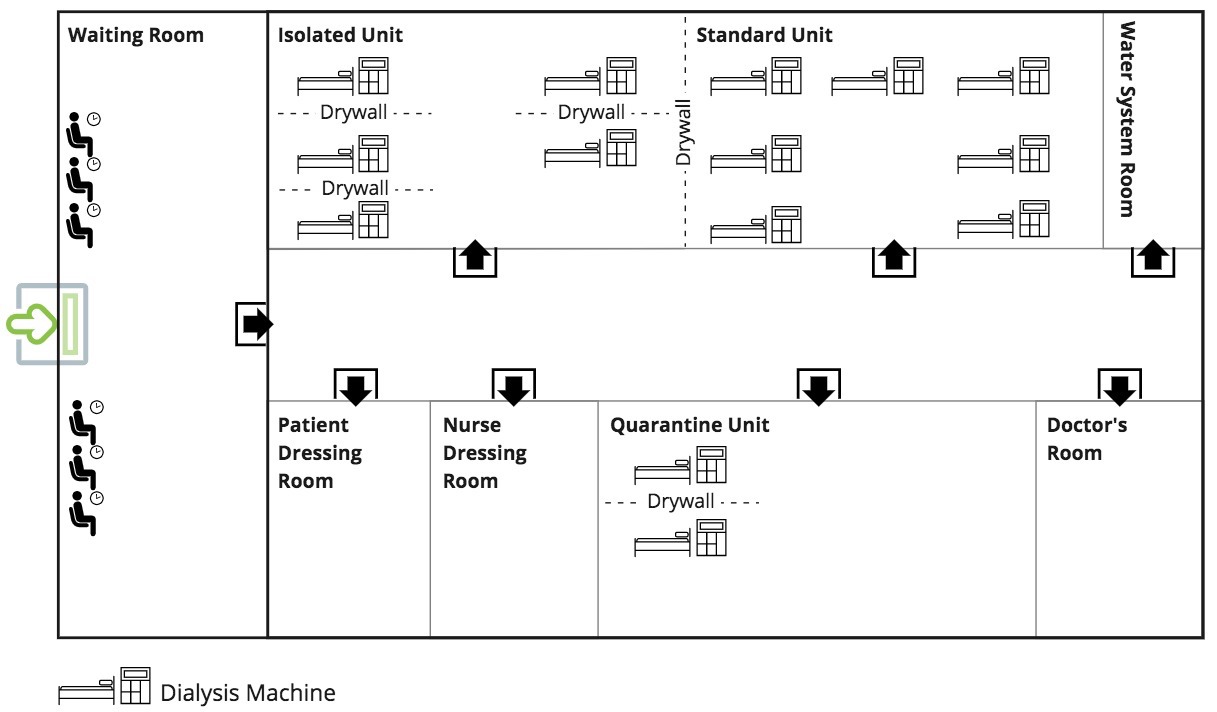} 
\label{FloorPlan}
\end{center}
\end{figure}

The dialysis clinic in our case hospital has 14 dialysis machines. Starting from 8:00 am, four sessions can be performed on each machine each day (Monday-Saturday). On Sundays, only emergency patients are treated in the clinic. Since the beginning of the pandemic, the hospital cohorts the hemodialysis patients into three cohorts (i.e., uninfected, infected, suspected), where each cohort is treated in a separate unit (standard, isolated and quarantine, respectively). Moreover, the units allocated to suspected and infected patients are divided into single-patient rooms, where each room includes one dialysis machine. All machines and rooms are cleaned intensively between dialysis sessions. To divide a unit into rooms, the hospital can build drywalls easily (overnight). A nurse can handle the treatment of four patients simultaneously in the same unit. 
Although the room arrangements within a unit can be made overnight, since all weekly plans (such as nurse shifts) must be fixed beforehand, the configuration of the units is arranged only at the end of the week using drywalls (i.e., Sunday). For instance, only two machines were allocated to infected patients during the first couple of months of the pandemics, while there exist five rooms in the infected patient unit (isolated unit) during the period we gathered the data from the hospital (November-December 2020). A schematic of the dialysis clinic during this period is given in Figure \ref{FloorPlan}. As shown in the figure, seven, five and two machines are currently allocated to standard, isolated and quarantine units, respectively. 

In addition to treating infected, suspected and uninfected patients in different units, the clinic schedules the treatments of different patient groups at different (nonoverlapping) sessions, as long as their capacities permit. Since all units are along the same corridor, patients and their companions use the same space before and after the dialysis sessions, and hence overlapping sessions may increase the risk of infection spread in the clinic. Therefore, the clinic arranges dialysis treatment schedules such that uninfected (chronic and acute) patients are assigned to morning sessions, which are followed by the sessions of suspected, and then the infected patients. Therefore, the number of rooms allocated to each unit at the beginning of a week directly affects the number of overlapping sessions that will be incurred each day. Currently, the hospital does not have any analytical tools that can support decisions for determining the size of units by considering daily treatment schedules. Given the uncertainties in the number of patients that will need dialysis treatments in a day, it is challenging to make weekly capacity planning decisions by considering their possible implications on daily treatment schedules and overlapping sessions. Our first objective in this study is to develop a mathematical model that can assist the hospital to determine cohort capacities effectively under demand uncertainty.

While our case hospital has separated the clinic into three units to treat three groups of patients, another cohorting strategy could be to divide the clinic into two units, where the noninfected patients are treated in the first unit, and the suspected and infected patients can be sequentially treated in the second unit, as practiced in some applications \citep{whiteside2020}. Depending on the number of patients over a week, which is not known beforehand (except the regular chronic patients), we investigate whether this alternative cohorting strategy might work better for our clinic to avoid overlapping sessions. Hence, our second objective in this paper is to evaluate the performance of different cohorting strategies and make recommendations to the hospital.

We next define the problem and then present the mathematical models developed to achieve the study's objectives.

\subsection {Problem Definition} \label{s:prob}
We consider the capacity planning problem 
of a dialysis clinic in a hospital setting, which allocates the valuable dialysis resources to serve different patient cohorts during a pandemic. Specifically, to mitigate infection spread among infected and uninfected patients, the hospital cohorts the patients and treat them in different units, and if capacity allows, in nonoverlapping sessions. There exist four types (groups) of patients with different characteristics, summarized in Table \ref{ptypes}.\\

\begin{table}[H] \setlength{\tabcolsep}{8pt}
\begin{center}
\begin{normalsize}
\begin{tabular}{ll}
\toprule
\textbf{Patient Types}&\textbf{Description}\\
\midrule
\textbf{Uninfected Acute} (Type 1) & Uninfected patients that are admitted to the hospital \\& and need dialysis. \\
\textbf{Uninfected Chronic} (Type 2) & Uninfected patients that receive dialysis treatment \\& regularly (every MWF or TTF).\\
\textbf{Infected COVID-19} (Type 3) & Infected patients with coronavirus that need \\& dialysis.\\
\textbf{Suspected COVID-19} (Type 4) & Suspected patients with the possibility of having \\& coronavirus that need dialysis.\\
\bottomrule
\end{tabular}
\end{normalsize}
\caption{Patient types treated in the dialysis clinic during the pandemic} 
\label{ptypes}
\end{center}
\end{table}

The daily demand for chronic dialysis patients (Type 2) can be considered fixed and known as these patients have to receive regular treatments. The number of uninfected acute patients (Type 1) that need dialysis each day is uncertain and can depend on a variety of factors such as the number of admitted patients in the hospital each day and their medical conditions which may cause temporary kidney failure. The daily demand by suspected and infected COVID-19 patients (Type 3 and 4) are also uncertain, which depend on uncontrollable factors including the current infection rate in the population, referrals from other clinics/hospitals, and the preferences of patients. Due to the high spread rate of COVID-19, a patient that was clean the day before can be suspected or infected the next day. Therefore, the daily patient status is important to consider in assigning patients to different units and setting their treatment times. In our model, we represent the uncertainties in the number of Type 1, 3 and 4 patients by a set of discrete scenarios. More specifically, each scenario specifies the estimated number of patients that will need dialysis each day over a week. In our case study (\S \ref{s:case}), we generate scenarios by making predictions from past patient demands.

Since COVID-19 has a high transmission risk, infected patients and uninfected patients should be treated in different units of the hospital. We evaluate two alternative cohorting strategies, which involve separating the clinic into three units and two units. We assume that the hospital management sets the cohorting strategy in advance as a policy. Table \ref{dutypes} shows the patient groups that are assigned to clinic units under each cohorting policy. It is not desirable for the dialysis clinic to turn down Type 3 and 4 patients as serving these patients is within the mandate of a pandemic hospital. Moreover, Type 1 patients are not usually transferable to other hospitals as they may be in critical condition and are already admitted to various hospital wards. Finally, Type 2 cannot be easily assigned to another hospital or clinic under the current regulations. Therefore, for a given cohorting policy, it is critical to make capacity configuration (planning) decisions in the clinic for the optimal use of valuable resources (i.e., dialysis machines) while providing care to all patients that need dialysis. We assume that a large penalty is incurred if a patient cannot be treated in the clinic due to capacity insufficiency. 

On each day, the clinic strives to provide treatment to different patient types at different periods to avoid physical interaction in the clinic as much as possible. When Type $i$ and Type $j$ patients from different cohorts are treated in the same session, we call this an ``overlap $i$x$j$". In particular, the clinic wants to avoid overlaps of uninfected patients with infected or suspected patients (overlaps 1x3, 1x4, 2x3, 2x4). Avoiding each of these overlaps is of primary importance for the clinic.
Moreover, in the three-unit cohorting strategy where the infected and suspected patients are treated in different rooms, their overlap (overlap 3x4) should be avoided if possible; however, this overlap is considered to be less serious compared to overlaps that involve clean patients. 
Note that in the two-unit cohorting strategy, Type 3 and 4 patients are sequentially treated in the same unit, and suspected patients are always treated earlier than infected patients. Since there is enough time in between sessions to disinfect the rooms, we assume there is no risk for a Type 3 patient to infect a Type 4 patient if they are consecutively scheduled in the isolated unit, hence this is not considered an overlap.

\begin{figure}[H]
\centering
\caption{Example daily treatment schedules under three-unit (a) and two-unit (b) cohorting. } \medskip
\includegraphics[width=\linewidth]{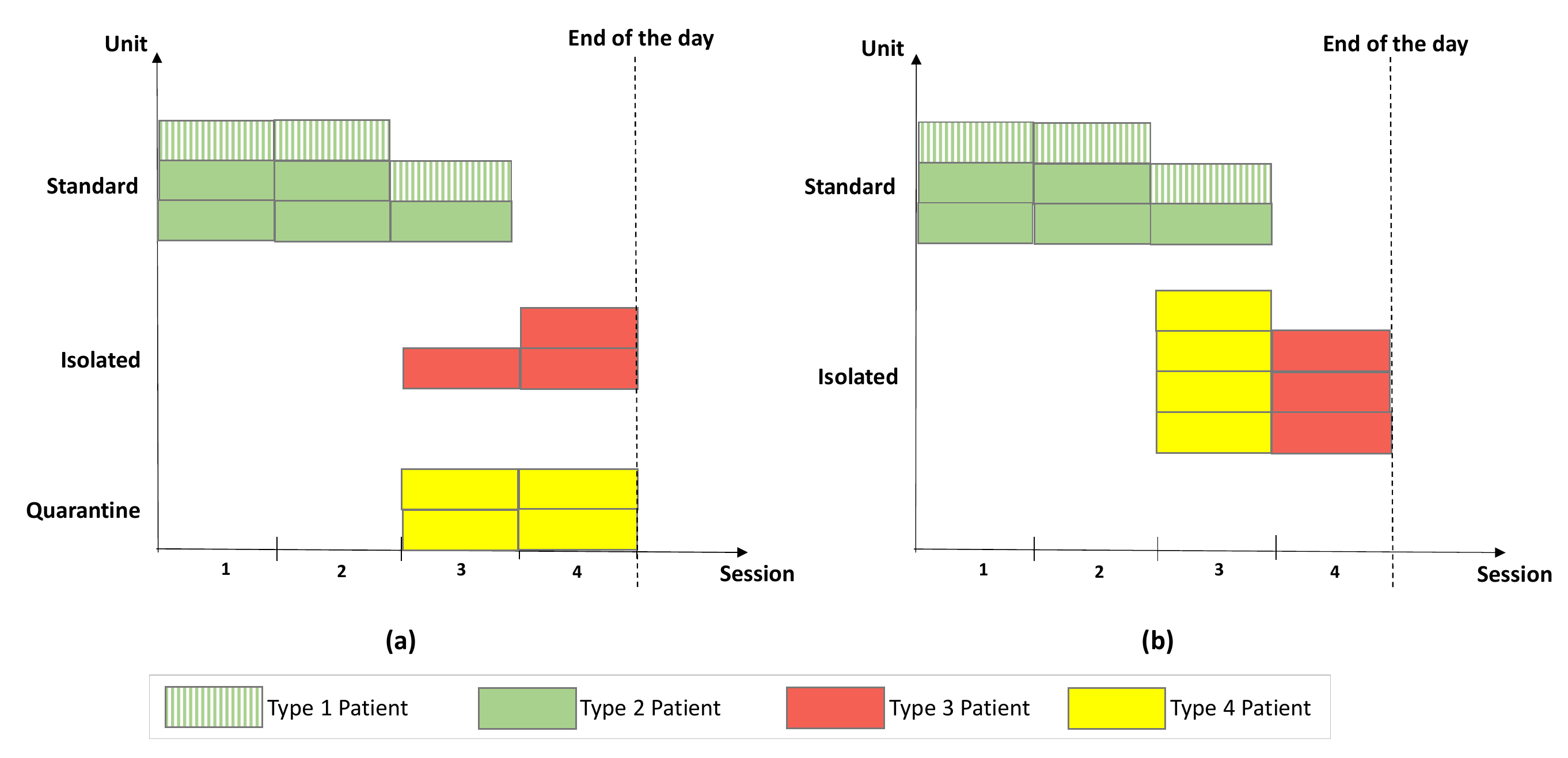}
\label{twothreerooms}
\end{figure}

Figure \ref{twothreerooms} illustrates example daily treatment schedules for both cohorting strategies on a small instance with seven machines, and fifteen patients. 
The number of patients of Type 1-4 are three, five, four and three, 
respectively. 
As shown in Figure \ref{twothreerooms}(a), when three, two, and two dialysis machines are allocated to standard, isolated and quarantine units, respectively, some patients from each type have to be served in overlapping sessions. In Figure \ref{twothreerooms}(b), three and four machines are allocated to standard and isolated units, respectively, and the sessions of some uninfected (Type 1 and 2) patients and suspected (Type 4) patients overlap. Therefore, both the cohorting policy and the capacity configuration decisions can have a significant effect on the number of patients who are exposed to infection risk each day.\\ 

\begin{table}[H] \setlength{\tabcolsep}{8pt}
\begin{center}
\begin{normalsize}
\begin{tabular}{lll}
\toprule
\textbf{Units}&\textbf{Three-Unit Cohorting}&\textbf{Two-Unit Cohorting}\\
\midrule
\textbf{Standard} & Uninfected Acute (Type 1) & Uninfected Acute (Type 1) \\& Uninfected Chronic (Type 2) & Uninfected Chronic (Type 2)\\
\hline
\textbf{Isolated} & Infected COVID-19 (Type 3) & Infected COVID-19 (Type 3)\\ & & Suspected COVID-19 (Type 4) \\
\hline \textbf{Quarantine} & Suspected COVID-19 (Type 4)\\
\bottomrule
\end{tabular}
\end{normalsize}
\caption{Three-unit and two-unit cohorting strategy patient assignments} 
\label{dutypes}
\end{center}
\end{table}

We consider a planning horizon of one week for the capacity planning decisions, which involve determining the number of dialysis machines in each unit. Capacity planning must be done at the beginning of each week, facing uncertain demand from each patient type over the following week. We follow a two-stage approach to formulate the clinic's problem, where capacity planning is done in the first stage, accounting for the scheduling decisions under different demand scenarios, given in the second stage. Accordingly, the capacity planning decisions are made at the beginning of each week, before uncertainties in demands are resolved by considering the implications of the second stage treatment scheduling decisions. Treatment scheduling decisions involve determining the timing of sessions allocated to each patient group each day given the predetermined first-stage decision variables and the realized patient demands. 

The objective is to minimize the total expected penalties that will occur over a week due to overlapping sessions. Specifically, a penalty is incurred for each patient that is scheduled in an overlapping session with other patients treated in different units. We impose an overlapping penalty charge per patient to account for the fact that there might be a higher infection risk if a larger number of conflicting patients are in the clinic at the same time. Additionally, if there are any patients that cannot be served due to capacity insufficiency, this causes infeasibility, which is modeled by imposing a large penalty cost. Different unit penalty charges might be set for not being able to treat different types of patients if it is more difficult to transfer certain types of patients to other facilities. 
Finally, to avoid alternative solutions that spread patients unnecessarily over multiple sessions, which may arise when capacity is abundant, we impose a small penalty for each session performed over a day.



The additional assumptions of our problem are as follows:
\vspace{-2mm}
\begin{itemize}
    \item Each dialysis unit can have a capacity restriction imposed by physical constraints. \vspace{-3mm}
    \item Only Type 1 and 2 patients can be served in the same unit during a session. Otherwise, a session must be composed of patients of the same type. \vspace{-3mm}
    \item The maximum number of sessions in each unit per day is equal. Moreover, each session has a fixed duration for all patient types. \vspace{-3mm}
    \item The time required for disinfecting and preparing the machines between sessions is fixed and included in a session's duration.\vspace{-3mm}
    \item One machine can serve one patient in each session.\vspace{-3mm}
    \item The sessions run in parallel; that is, the starting and ending times of all sessions are the same.\vspace{-3mm}
    \item In each unit, the dialysis sessions are scheduled consecutively to allow easier planning of other resources (e.g., nurses). That is, there is a no-idling policy. \vspace{-3mm}
    \item The sessions of uninfected patients are scheduled at earlier sessions starting from the beginning of the day, and if capacity allows, the sessions of infected and suspected patients start after all infected patients are treated.\vspace{-3mm}
    \item In the two-unit cohorting, the sessions of suspected patients precede those of infected patients in the isolated unit.\vspace{-3mm}
\end{itemize}
We next present stochastic programming models to solve the capacity planning problem of the hospital described here.

\subsection{Mathematical Models} \label{s:model}

We present two integer programming models for the alternative cohorting strategies that separate the dialysis clinic into units. We first present the notation and then present the models for three- and two-unit models in \S 3.3.1 and \S 3.3.1, respectively. \\

\vspace{-2mm}
\noindent \textbf{Sets} \\
$K$: set of scenarios; $k \in K$ \\
$I$: set of patient types; $i \in I$; $i$=1,..,4 \\
$J$: set of dialysis units; $j \in J$; $j$=1,..,3 (1: Standard unit, 2: Isolated unit, 3: Quarantine unit) \\
$T$: set of days in the planning horizon; $t \in T$; $t$=1,..,6 \\
$S$: set of dialysis sessions on each day; $s \in S$; $s$=1,..,4 \\

\vspace{-2mm}
\noindent \textbf{Parameters} \\
$p_{k}$: probability of scenario k\\
$H_{it}^k$: number of type $i$ patients who need to receive dialysis on day $t$ of scenario $k$ \\
$C_j$: maximum number of dialysis machines that can be allocated to unit $j$\\
$\hat{C}$: total number of available dialysis machines  \\
$\alpha_1$: penalty (per patient) for uninfected patients and infected COVID-19 patients that are treated in overlapping sessions \\
$\alpha_2$: penalty (per patient) for uninfected (Type 1, 2) patients and suspected COVID-19 (Type 4) patients that are treated in overlapping sessions \\
$\alpha_3$: penalty (per patient) for the conflict of suspected (Type 4) and infected (Type 3) patients that are treated in overlapping sessions \\
$\Pi_i$: penalty (per patient) for infeasibility if a Type $i$ patient cannot be treated due to capacity insufficiency. \\
$\epsilon$: small penalty coefficient for starting a dialysis session (to avoid unnecessary sessions) \\

\vspace{-2mm}
\noindent \textbf{First stage decision variables} \\
$R_j$: number of dialysis machines allocated to dialysis unit $j$\\

\vspace{-2mm}
\noindent \textbf{Second stage decision variables} \\
$X_{its}^k$: number of Type $i$ patients scheduled to receive dialysis in session $s$ on day $t$ of scenario $k$ \\
$F_{it}^k$: number of Type $i$ patients that cannot be served on day $t$ of scenario $k$ \\
$N_{jts}^k$:	$\begin{cases}
  1, & \text{if dialysis unit type $j$ is used for session $s$ on day $t$ of scenario $k$,} \\
  0, & \text{otherwise}.
\end{cases}$ \\
$Q_{ts}^k$: number of uninfected (Type 1, 2) and infected (Type 3) patients treated in session $s$ on day $t$ of scenario $k$ \\
$G_{ts}^k$: number of uninfected (Type 1, 2) and suspected (Type 4) patients treated in session $s$ on day $t$ of scenario $k$ \\
$W_{ts}^k$: number of infected (Type 3) and suspected (Type 4) patients treated in session $s$ on day $t$ of scenario $k$ \\ 
$U_{ts}^k$: $\begin{cases}
  1, & \parbox[t]{0.9\textwidth}{if an overlap exists in session $s$ on day $t$ of scenario $k$ for uninfected (Type 1, 2) and infected (Type 3) patients,} \\
  0, & \text{otherwise}.
\end{cases}$ \\
\noindent $D_{ts}^k$:	$\begin{cases}
  1, & \parbox[t]{0.9\textwidth}{if an overlap exists in session $s$ on day $t$ of scenario $k$ for uninfected (Type 1, 2) and suspected (Type 4) patients} \\
  0, & \text{otherwise}.
\end{cases}$ \\
$V_{ts}^k$:	$\begin{cases}
  1, & \parbox[t]{0.9\textwidth}{if an overlap exists in session $s$ on day $t$ of scenario $k$ for infected (Type 3) and suspected (Type 4) patients,} \\
  0, & \text{otherwise}.
\end{cases}$ 
$Z_{3u}$: auxiliary variable for the objective function 



\subsubsection{Three-unit cohorting model} \label{3unitmodel}
The formulation for making capacity allocation decisions to implement the three-unit cohorting strategy is as follows.   
\begin{equation}
\text{min} \quad  Z_{3u} = \sum_{k \in K} p_k \bigg[ \alpha_1 \sum_{t\in T, s\in S}  Q_{ts}^k + \alpha_2 \sum_{t\in T, s\in S}  G_{ts}^k + \alpha_3 \sum_{t\in T, s\in S}  W_{ts}^k 
+ \Pi_i \sum_{i\in I, t\in T}  F_{it}^k + \epsilon \sum_{j\in J, t\in T, s\in S}N_{jts}^k \bigg] \label{obj1} 
\end{equation}
\qquad \qquad subject to \vspace{-4mm}
\begin{align}
R_j \leq C_j \qquad \qquad \qquad &\forall j \in J \label{c9}\\
\sum_j R_j \leq \hat{C} \qquad \qquad \qquad & \phantom{x} \label{c10}\\
X_{1ts}^k + X_{2ts}^k \leq R_1 \qquad \qquad \qquad &\forall t \in T, s \in S, k \in K \label{c3}\\
X_{3ts}^k \leq R_2 \qquad \qquad \qquad &\forall t \in T, s \in S, k \in K \label{c4}\\
X_{4ts}^k \leq R_3 \qquad \qquad \qquad &\forall t \in T, s \in S, k \in K \label{c5}\\
\sum_{s} X_{its}^k = H_{it}^k -F_{it}^k \qquad \qquad \qquad &\forall i \in I, t \in T, k \in K \label{c1}\\
X_{1ts}^k + X_{2ts}^k \leq M\times N_{1ts}^k \qquad \qquad \qquad &\forall t \in T, s \in S, k \in K \label{c6}\\
N_{1ts}^k \geq N_{1t(s+1)}^k \qquad \qquad \qquad &\forall s \in S \colon s\leq3,  t \in T,  k \in K \label{c11}\\
X_{3ts}^k \leq M\times N_{2ts}^k \qquad \qquad \qquad &\forall t \in T, s \in S, k \in K \label{c7}\\
X_{4ts}^k \leq M\times N_{3ts}^k \qquad \qquad \qquad &\forall t \in T, s \in S, k \in K \label{c8}\\
\sum_{s'=s+2}^{S} N_{jts'}^k \leq |S|(1-N_{jts}^k+N_{jt(s+1)}^k)\qquad \qquad \qquad &\forall s \in S \colon s\leq2,  t \in T,   j \in J \colon j\neq1,  k \in K \label{c12} \\
1+U_{ts}^k \geq N_{1ts}^k+N_{2ts}^k \qquad \qquad \qquad &\forall t \in T, s \in S, k \in K \label{c13}\\
Q_{ts}^k \geq X_{1ts}^k+X_{2ts}^k + X_{3ts}^k- M(1-U_{ts}^k) \qquad \qquad \qquad &\forall t \in T, s \in S, k \in K \label{c14}\\
1+V_{ts}^k \geq N_{2ts}^k+N_{3ts}^k \qquad \qquad \qquad &\forall t \in T, s \in S, k \in K \label{c15} \\
W_{ts}^k \geq X_{3ts}^k + X_{4ts}^k - M(1-V_{ts}^k) \qquad \qquad \qquad &\forall t \in T, s \in S, k \in K \label{c16}\\
1+D_{ts}^k \geq N_{1ts}^k+N_{3ts}^k \qquad \qquad \qquad &\forall t \in T, s \in S, k \in K \label{c17}\\
G_{ts}^k \geq X_{1ts}^k + X_{2ts}^k  + X_{4ts}^k- M(1-D_{ts}^k) \qquad \qquad \qquad &\forall t \in T, s \in S, k \in K \label{c18}\\
R_{j} \in \mathbb{Z}^+  \qquad \qquad \qquad &\forall j\in J \label{c23}\\
X_{its}^k \in \mathbb{Z}^+  \qquad \qquad \qquad &\forall i \in I,  t \in T, s \in S, k \in K \label{c20}\\
F_{it}^k \in \mathbb{Z}^+  \qquad \qquad \qquad &\forall i \in I, t \in T,  k \in K \label{c21}\\
Q_{ts}^k, W_{ts}^k, G_{ts}^k \in \mathbb{Z}^+  \qquad \qquad \qquad &\forall t \in T, s \in S, k \in K \label{c22}\\
N_{jts}^k \in \{0, 1\}  \qquad \qquad \qquad &\forall j \in J,  t \in T, s \in S, k \in K \label{c24}\\
U_{ts}^k, V_{ts}^k, D_{ts}^k \in \{0, 1\}  \qquad \qquad \qquad &\forall t \in T, s \in S, k \in K \label{c25}
\end{align}

The first three terms of the objective function \eqref{obj1} minimizes the weighted sum of the penalties due to overlapping sessions of patient groups that are not preferred to be in the clinic at the same times to mitigate the infections transmission risk. The fourth and fifth terms in \eqref{obj1} are for penalizing patients that cannot be treated in the clinic due to capacity insufficiency. The last term in the objective function minimizes the number of dialysis sessions provided during a week to ensure that the sessions are scheduled in a compact way. Constraints \eqref{c9} ensure that the number of machines allocated to each unit does not exceed its  capacity. Constraint \eqref{c10} ensures that the total number of dialysis machines assigned to units does not exceed the number of machines available in the clinic. Constraints \eqref{c3}-\eqref{c5} guarantee that for any unit, the number of patients assigned to a session cannot exceed the capacity reserved for that unit. Constraints \eqref{c1} 
balance the daily number of dialysis sessions required and provided. Constraints \eqref{c6} and \eqref{c11} ensure that dialysis sessions of uninfected (Type 1 and 2) patients start from the beginning of the day, and multiple sessions are conducted consecutively in the standard unit. Constraints \eqref{c7}, \eqref{c8} are for determining used sessions of suspected (Type 4) and infected (Type 3) patients anytime during a day; additionally, constraints \eqref{c12} ensure that the sessions in their corresponding units must be conducted consecutively for the compactness of the schedule. In other words, by keeping track of the two consecutive sessions, \eqref{c12} prevents leaving an idle session between any sessions in which suspected (Type 4) and infected (Type 3) patients are treated. Constraints \eqref{c13} determine the overlapping sessions among standard and isolated units. Accordingly, the number of uninfected (Type 1 and 2) and infected (Type 3) patients that receive dialysis in the same session is determined by constraints \eqref{c14}. Similarly, the number of infected (Type 3) and suspected (Type 4) patients that receive dialysis in the same session is determined by constraints \eqref{c15} and \eqref{c16}; further, the number of uninfected (Type 1 and 2) and suspected (Type 4) patients that receive dialysis in the same session is determined by constraints \eqref{c17} and \eqref{c18}. Finally, the integer variables are defined by constraints \eqref{c23}-\eqref{c22}, and binary variables are defined by constraints \eqref{c24} and \eqref{c25}.

In our computations, we set the value of $M$ in constraints \eqref{c6}, \eqref{c7}, and \eqref{c8} equal to the number of dialysis machines available, $\hat{C}$. The $M$ values in constraints \eqref{c14}, \eqref{c16}, and \eqref{c18} are set as the number of dialysis treatments that can be provided by the clinic, which equal to the multiplication of the number of machines available and the number of sessions on a day ($|S| \times \hat{C}$).

\subsubsection{Two-unit cohorting model} \label{2unitmodel}
The model for the two-unit cohorting strategy has a few differences compared to the model presented above for the three-unit strategy. First, there exist two units, that is, $j=1,2 \in J$, where $j=1$ denotes the standard unit where the uninfected patients are treated, and $j=2$ denotes the isolated unit where suspected and infected patients are treated sequentially. To sequence the treatments of infected patients after the suspected patients in the isolated unit, we introduce a new variable $Y_{ts}^k$, which takes the value of 1 when all suspected patients that can be treated in a day are assigned to a dialysis session, thereby specifying the earliest session that the infected patients can be scheduled in the unit. Moreover, a dummy session is needed to schedule consecutive sessions of suspected and infected patients in the same cohort by using variables $Y_{ts}^k$. We denote the dummy session by $\{0\}$ and let $S_0 = S \cup \{0\}$. The $W_{ts}^k$ and $V_{ts}^k$ variables that are used before in the three-unit cohorting model are not relevant in the two-unit model, since the suspected and infected patients receive in the same unit in separate sessions. Finally, we define $Z_{2u}$ as an auxiliary variable to keep the objective function value attained for the two-unit model.

The model for the capacity planning problem of the clinic to implement a 2-unit cohorting strategy is presented below.
\begin{equation}
\text{min} \qquad  Z_{2u} =\sum_{k\in K} p_k \bigg[ \alpha_1 \sum_{t\in T, s \in S}  Q_{ts}^k + \alpha_2 \sum_{t \in T, s\in S}  G_{ts}^k + 
\Pi_i \sum_{i\in I, t\in T}  F_{it}^k + \epsilon \sum_{j\in J, t\in T, s\in S}N_{jts}^k \bigg] \label{obj2}
\end{equation}
\qquad \qquad \qquad \qquad subject to \vspace{-2mm}
\begin{align}
\eqref{c9}, \eqref{c10}, \eqref{c3}, \eqref{c1}, 
\eqref{c6}, \eqref{c11}, \eqref{c13}, \eqref{c14} \qquad & \phantom{x} \notag\\
\eqref{c17}, \eqref{c18}, \eqref{c23}, \eqref{c20}, \eqref{c21}, \eqref{c24} \qquad & \phantom{x} \notag\\
X_{3ts}^k + X_{4ts}^k \leq R_2 \qquad \qquad \qquad &\forall t \in T, s \in S, k \in K \label{c29}\\
X_{3ts}^k + X_{4ts}^k \leq M \times N_{2ts}^k \qquad \qquad \qquad &\forall t \in T, s \in S, k \in K \label{c31} \\
\sum_{s'=s+2}^{4} N_{2ts'}^k \leq |S|(1-N_{2ts}^k+N_{2t(s+1)}^k)\qquad \qquad \qquad &\forall s \in S \colon  s\leq2,  t\in T,  k \in K \label{c35}\\
H_{4t}^k - \sum_{s'=0}^s X_{4ts'}^k - F_{4t}^k \leq M\times Y_{ts}^k \qquad \qquad \qquad &\forall t \in T,  s\in S_0,  k \in K \label{c40}\\
X_{3ts}^k \leq M\times (1-Y_{t(s-1)}^k) \qquad \qquad \qquad &\forall t \in T, s \in S, k \in K \label{c41} \\
X_{it0}^k = 0 \qquad \qquad \qquad &\forall i \in I,  t \in T,  k \in K \label{c42}\\
Q_{ts}^k, G_{ts}^k \in \mathbb{Z}^+  \qquad \qquad \qquad &\forall t \in T, s \in S, k \in K \label{c45} \\
U_{ts}^k, D_{ts}^k, Y_{ts}^k \in \{0, 1\}  \qquad \qquad \qquad &\forall t \in T, s \in S, k \in K \label{c48}
\end{align}

The objective function \eqref{obj2} is similar to \eqref{obj1}, except that there is no penalty now associated with overlapping Type 3 and 4 patients, which are sequentially treated in the same unit in the two-unit cohorting strategy. Moreover, constraints \eqref{c4} and \eqref{c5} in the three-unit model are consolidated in a single constraint \eqref{c29} in the two-unit model since Type 3 and 4 patients are assigned to the same unit. Similar to constraints \eqref{c7} and \eqref{c8}, constraints \eqref{c35} ensure that dialysis sessions for Type 3 and 4 infected patients can start anytime during a day, but their sessions must be conducted consecutively. Constraints \eqref{c40} and \eqref{c41} specify the earliest session for Type 3 patients after Type 4 patients are scheduled. Constraints \eqref{c48} define the domains of variables. The remaining constraints are identical to those of the three-unit model.

In implementing the two-unit model, the value of $M$ in constraints \eqref{c31} is set equal to $\hat{C}$. In constraints \eqref{c40} and \eqref{c41}, $M$ values are set to the largest number of Type 4 and Type 3 demand occurrences among $|K|$ scenarios, respectively.

\section{Case Study: Analysis and Results} \label{s:case}
In this section, we present a case study to test and illustrate the proposed models based on data obtained from the dialysis clinic in our collaborating hospital, whose operations are described before (\S \ref{s:system}). In \S \ref{s:case2}, we present an analysis of the data obtained from the dialysis clinic and explain how we generate demand scenarios to implement the proposed models. In \S \ref{s:case3}, we present a computational study and discuss the practical implications of the results.

\subsection{Data Analysis and Case Instances} \label{s:case2}
The hemodialysis clinic of our collaborating hospital has been applying a three-unit cohorting strategy since the beginning of the pandemic. The clinic does not have an advanced information system that keeps data regarding patient demands. Upon our request, a nurse working in the hemodialysis clinic of the Marmara University Pendik Hospital has recorded patient data for eight weeks in November and December 2020, which corresponds to a period around the second peak of the pandemic. The acquired data from the clinic include the number of patients of each type treated in the clinic each day. 

\begin{figure}[!hbt] 
\begin{center}
\caption{Daily number of dialysis patients by type}
\medskip
\includegraphics[scale=0.50]{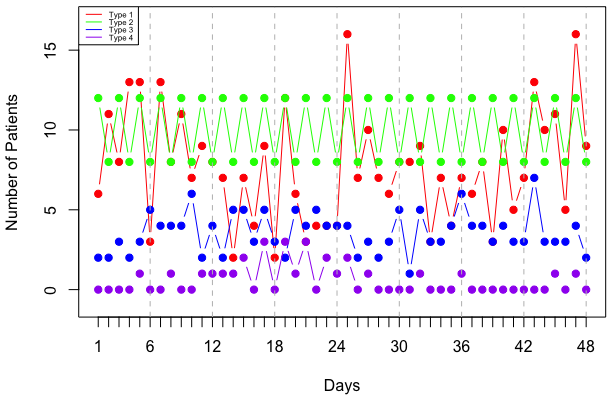} 
\label{graph:dailynopatient}
\end{center}
\end{figure}

We consider the eight weeks of data collected between November 2 to December 26 by excluding Sundays, and obtain a data set for 48 days in total. The number of patients is very low on Sundays compared to other days because only extremely urgent patients of Type 1 or Type 3 are treated on those days, and unlike rest of the days with five nurses working in the clinic, only one nurse is working on Sundays. The average number of patients being treated is four, on Sundays, which is extremely low compared to the weekdays and Saturdays. Therefore, we treat Sunday as an outlier and remove it from the whole data set. Figure \ref{graph:dailynopatient} presents a plot of the daily number of patients for each type and weekly intervals are indicated by the grey dashed lines. Besides, the raw data and descriptive statistics are presented in Appendix A. We observe that the daily number of uninfected acute patients (Type 1) is highly uncertain, as shown by fluctuations in the plot. On the other hand, the number of uninfected chronic patients (Type 2) shows a periodic pattern and constant throughout the entire time horizon. The infected COVID-19 (Type 3) and suspected COVID-19 (Type 4) patients are a smaller percentage of the overall demand, while the number of infected COVID-19 (Type 3) patients also varies considerably over the observed period.

We create demand scenarios based on the available data to illustrate our models.
Since the unit configuration can be changed once per week, we consider a weekly (six days) planning horizon. Therefore, given the realized demands in the past, demand scenarios should be generated for the following week. One plausible approach to estimate demand coming from infected (Type 3) or suspected (Type 4) COVID-19 patients is to consider the dynamics of the pandemic in the country. SIR-based disease progression dynamic models have been used to estimate the demand for hospital services (cf. \cite{weissman2020locally}). An approach like this could be to use the aggregate prevalence estimates from such dynamic models and assume a fraction of those patients will need dialysis. However, while these models work well at the aggregate level, they may not reliably estimate the patients of a single dialysis clinic, given that there exist other clinics serving the population. Therefore, for predicting the demand for each type of patient we use different methods according to the characteristics of the available data.\\

\begin{figure}[!hbt]
\begin{center}
\caption{Total number of Type 3 and 4 patients treated in the hemodialysis clinic and the reported national COVID-19 cases (November-December 2020)} 
\medskip
\label{fig:CovidTrendImage}
\includegraphics[scale=0.35]{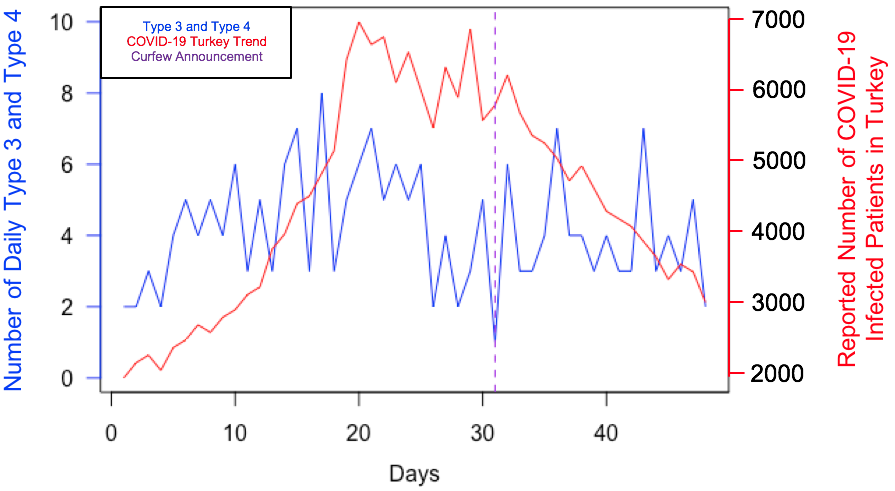} 
\end{center}
\end{figure}

Figure \ref{fig:CovidTrendImage} shows the number of infected (Type 3) and suspected (Type 4) COVID-19 patients treated in our case hospital and the daily reported COVID-19 infected patients in Turkey between November 2 and December 26 (excluding Sundays). We observe that the number of Type 3 and 4 patients treated in the pandemic hospital does not closely follow the same pattern with the number of cases in the country. When the COVID-19 cases peaked in the middle of November, the numbers of hemodialysis patients with COVID-19 increased and peaked as well. However, after a curfew was announced at the beginning of December, which is indicated with the purple dashed line in the figure, the number of nationwide cases started to diminish, while the number of hemodialysis patients with COVID-19 seems stationary. There is also a  higher variability in the demand for dialysis, which is at a smaller scale compared to the aggregate number of cases at the national level. 
Time series methods are expected to work well for short-term forecasting during the pandemic \citep{doornik2020short}. Therefore, we use time series forecasting methods to estimate the demand for the upcoming week. 


\newpage
Based on the clinic data, we generate scenarios that represent weekly patient demand estimates. To obtain weekly demand estimates, we divide the data set with 48 observations (days) into two, namely, a training set and a testing set for validation purposes. The last week (six days) of the data points are used as a testing set whereas the first 42 data points are used as the training data set and the forecasting methods are applied to first seven weeks of the data. The accuracy of the estimators is verified by comparing the estimators with the testing set that contains six data points. Furthermore, this process is repeated by producing estimations for week six and seven based on the data that corresponds to the first five and six weeks, respectively. We next explain the procedure followed to predict demands for each patient type.

\begin{figure}[!hbt]
\caption{Prediction intervals for Type 1 (a), Type 3 (b) and Type 4 (c) patients for week 8} \medskip \medskip
\label{fig:abc}
\begin{minipage}[c]{0.48\linewidth} 
\centering 
\includegraphics[width=7.5cm]{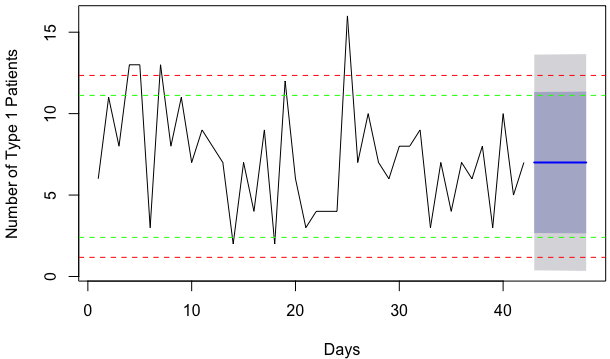} 
(a)
\label{T1SES}
\end{minipage}
\hfill
\begin{minipage}[c]{0.48\linewidth} 
\centering
\includegraphics[width=7.5cm]{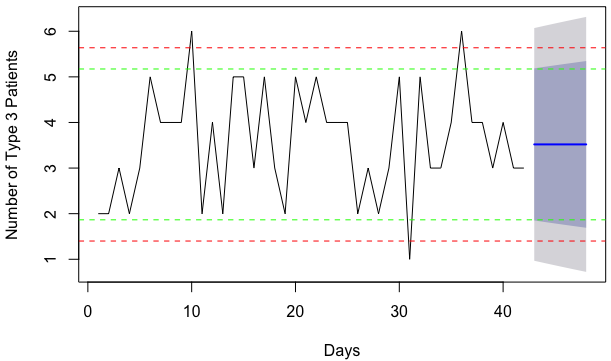} 
(b)
\label{T3SES}
\end{minipage}
\end{figure}

\vspace{-2mm}
\begin{itemize}
\item \textit{Uninfected Chronic Patients (Type 2)}. According to the data, the number of uninfected chronic patients is constant through the entire horizon, 12 patients are treated according to the Monday-Wednesday-Friday (MWF) regime and eight patients are enrolled to Tuesday-Thursday-Saturday (TTS) regime. We set the future daily demands for Type 2 patients by considering these actual number of chronic patients regularly treated by the clinic.

\item \textit{Uninfected Acute (Type 1), Infected COVID-19 (Type 3) and Suspected COVID-19 (Type 4) Patients}. To predict the demand coming from these patient types, we use simple exponential smoothing method to minimize the Root Mean Square Errors (RMSE) and find the point estimate for the next week's daily demand for patient type $i$, $\hat{Y_i}$. 
The $80\%$ and $90\%$ prediction intervals are calculated as $\hat{Y_i}\pm 1.28\sigma_i$ and $\hat{Y_i}\pm 1.64\sigma_i$ respectively, where the standard deviation $\sigma_i$ is approximated by the RMSE. The RMSE values are multiplied by constant $k$, which are 1.28 and 1.64, in order to create $80\%$ and $90\%$ prediction intervals (PI), respectively, to forecast Type $i$ demands. Finally, assuming a uniform distribution with each prediction interval, we discretize the distributions by rounding negative values to 0 and the real numbers to the nearest integer. This process is applied to first five, six and seven weeks to create discrete distributions for the following week (i.e., next six days). The estimated distributions are used to create instances to provide a basis for possible scenarios that might occur in the dialysis clinic. The 80\% and 90\% prediction intervals can be considered representing a less and more conservative approach in accounting for uncertainty, respectively. The prediction intervals for each patient type for each week are presented in Table \ref{table:bounds}. 

Figure \ref{fig:abc} depicts the data and the prediction intervals for Type 1, Type 3 and Type 4 patients respectively, using seven weeks of past data to predict week eight. The blue line represents the point estimate, green dashed lines represent $80\%$ PI approximation 
and red dashed lines represent $90\%$ PI approximation. 
Discrete distributions derived for the demands of Type 3 and Type 4 patients in week 8 are in Appendix B as an example. \\
\end{itemize}

\begin{table}[hbt!]
\small
\centering
\begin{tabular}{|l|ll|ll|ll|}
\hline
\textbf{} & \multicolumn{2}{c}{\textbf{Week 6}} & \multicolumn{2}{c}{\textbf{Week 7}} & \multicolumn{2}{c|}{\textbf{Week 8}} \\
\hline
\textbf{Patient Type} & \textbf{80\% PI}  & \textbf{90\% PI}  & \textbf{80\% PI}  & \textbf{90\% PI}   & \textbf{80\% PI}  & \textbf{90\% PI}   \\ \hline
Type 1 & (2.31, 11.94) & (0.96, 13.29) & (2.30, 11.40) & (1.01, 12.68) & (2.40, 11.12) & (1.18, 12.65) \\
Type 3 & (1.87, 5.18) & (1.41, 5.64) & (2.20, 5.71) & (1.70, 6.21) & (1.87, 5.17) & (1.40, 5.64) \\
Type 4 & (-0.73, 1.77) & (-1.08, 2.21) & (-0.76, 1.59) & (-1.09, 1.92) & (-1.01, 1.17) & (-1.32, 1.48) \\
\hline
\end{tabular}
\caption{Demand prediction intervals for weeks 6, 7 and 8}
\label{table:bounds}
\end{table}

\begin{figure}[!hbt] \begin{minipage}[c]{0.48\linewidth}
\centering
\includegraphics[width=7.5cm]{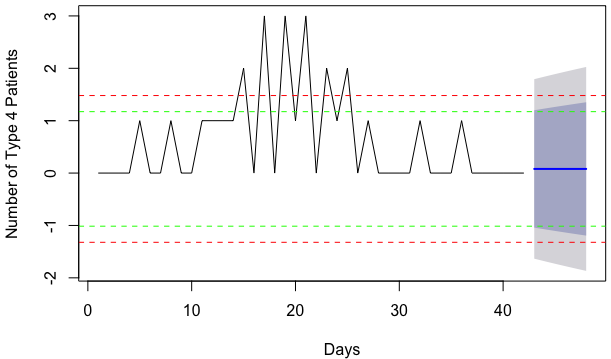} 
(c)
\label{T4SES}
\end{minipage}
\hfill
\begin{minipage}[c]{0.48\linewidth} 
\centering
\includegraphics[width=4.5cm]{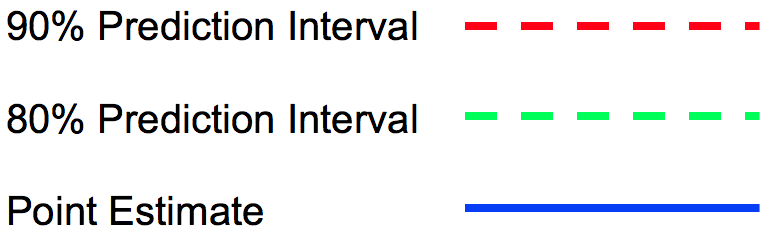} 
\label{T3SES}
\end{minipage}
\end{figure}

We generate 30 equiprobable scenarios to test our models. Demands per day in each scenario are generated by sampling from the corresponding probability distributions for each patient type (see Appendix B for the probability distributions). The other parameters used in our case instances are as follows. The total number of machines, $\hat{C}$ is set at 14. The capacities of units are set as $R_1$=11, $R_2$=8 and $R_3$=5 for the three-unit cohorting model, and as $R_1$=11 and $R_2$=8 for the two-unit cohorting model. In both models, the objective function weights are set as $\alpha_1$=$\alpha_2$=1,000, $\alpha_3$=100, $\Pi_i$=100,000 and $\epsilon$=2.

The proposed models are solved on a computer with an Intel(R) Core (TM) i7-9750H CPU @ 2.60 GHz processor and a 16 GB RAM. Gurobi 9.0.3 is used to solve the instances. We next present the results.

\subsection{Results and Analysis} \label{s:case3}

In this section, we present our results and analysis based on the case data. In \S \ref{s:case3a}, we evaluate the effectiveness of the hospital's current capacity configuration. In \S \ref{s:case3b}, we evaluate the solutions of the proposed stochastic optimization model for the three-unit cohorting strategy. In \S \ref{s:case3c}, we compare the solutions and performance of the alternate cohorting policies.


\subsubsection{Evaluation of the current capacity configuration} \label{s:case3a}

As described in \S \ref{s:system}, our case hospital applies a three-unit cohorting strategy to minimize the spread risk of COVID-19 among its patients. During November-December 2020, the clinic allocated seven dialysis machines to the standard unit to treat uninfected (Type 1 and 2) patients, two machines to the quarantine unit to treat suspected COVID-19 (Type 4) patients, and five machines to the isolated unit to serve the infected COVID-19 (Type 3) patients.

We first determine the number of weekly overlaps that would occur under the hospital's predetermined capacity configuration. For this purpose, for each week, we solve the deterministic version (i.e., with perfect demand information) of our three-unit cohorting model, \eqref{obj1}-\eqref{c25}, by fixing the number of resources allocated to the three units and using the realized weekly demands (given in Appendix A) as a single demand scenario. 
Then we find the optimal capacity allocation by solving the two-stage optimization model with this one scenario. Comparing the two solutions will allow us to illustrate the impact of capacity allocation decisions on the overlaps. We expect an improvement over the hospital's decision due to the assumption of perfect information and optimization of the capacity allocation. Table \ref{tablev1} compares the number of overlapping patients under the current and optimal capacity allocation decisions for each type of overlap $ixj$, which is denoted by $O_{ixj}$. For example, $O_{(1,2)x3}$ indicates the number of overlaps among uninfected (Type 1 and 2) patients with infected (Type 3) patients. 

\begin{table}[hbtp] \setlength{\tabcolsep}{4pt}
\scriptsize
\begin{center}
\begin{footnotesize}
\begin{tabular}{c|cccc|ccc|cccc|c}
\cmidrule[\heavyrulewidth]{2-12}
\multicolumn{1}{c}{} & \multicolumn{4}{l|}{\textbf{Performance with}} & \multicolumn{3}{l}{\textbf{Optimal}}&\multicolumn{5}{|l}{\textbf{Performance with optimal}} \\
\multicolumn{1}{c}{} & \multicolumn{4}{l|}{\textbf{hospital's allocation}} & \multicolumn{3}{l|}{\textbf{allocation }} &\multicolumn{5}{l}{\textbf{allocation (known demand)}} \\
\toprule
Week & $O_{(1,2)x3}$ & $O_{(1,2)x4}$ & $O_{3x4}$ & $Z_{3u}$ & $R_1$ & $R_2$ & $R_3$ & $O_{(1,2)x3}$ & $O_{(1,2)x4}$ & $O_{3x4}$ & $Z'_{3u}$ & $\frac{Z_{3u}-Z'_{3u}}{Z_{3u}}$  \\
\midrule
1 & 7 & 5 & 4 & 12,450 &     10 & 3 & 1 & 0 & 0 & 4 & 444 &     96\% \\
2 & 16 & 0 & 13 & 17,360 &   9 & 4 & 1 & 0 & 0 & 3 & 350 &      98\% \\
3 & 0 & 8 & 17 & 9,752 &     7 & 5 & 2 & 0 & 8 & 17 & 9,752 &   0\%  \\
4 & 5 & 15 & 15 & 21,558 &   8 & 3 & 3 & 0 & 0 & 17 & 1,758 &   92\% \\
5 & 15 & 11 & 4 & 26,456 &   8 & 4 & 2 & 8 & 6 & 10 & 15,050 &  43\% \\
6 & 2 & 0 & 12 & 3,254 &     9 & 4 & 1 & 0 & 0 & 5 & 546 &      83\% \\
7 & 0 & 0 & 0 & 48 &         10 & 4 & 0 & 0 & 0 & 0 & 36 &      25\% \\
8 & 34 & 11 & 4 & 45,458 &   9 & 4 & 1 & 15 & 2 & 19 & 17,950 & 61\% \\
\midrule
Total & 79 & 50 & 69 & 136,336 & - & - & - & 23 & 16 & 75 & 45,886 & 66\% \\
\bottomrule
\end{tabular}
\end{footnotesize}
\end{center} 
\caption{Comparison of the overlaps under hospital's current capacity allocation ($R_1=7, R_2=5, R_3=2$) with optimal allocation under perfect demand information.}
\label{tablev1}
\end{table}




As shown in Table \ref{tablev1}, except for week 3, there is a difference in resource allocation decisions, which indicates that dialysis machines could be better allocated in the past, if the demand was better estimated and the allocation was optimized. 
Results show that allocation decisions can significantly impact the overlaps of uninfected and infected patients. We discuss the mechanisms in the following. For all weeks except week 3, the hospital reserves fewer dialysis machines for the uninfected patients than they should (i.e., optimal $R_1$ is greater than 7).  As a result, since capacity at earlier sessions is not enough, some of the uninfected patients receive dialysis at later sessions of the day, during which suspected and infected patients also receive their treatments. 
As a result, under the hospital's allocation policy, the number of overlaps in uninfected and suspected patients' treatments increase, which results in greater objective values. The hospital's allocation is optimal only in week 3, in which demands from uninfected acute patients and suspected patients are lower than average, whereas a larger number of sessions is needed for infected patients. That is, high demands from infected patients in the third week appear to justify the capacity allocation policy of the hospital.

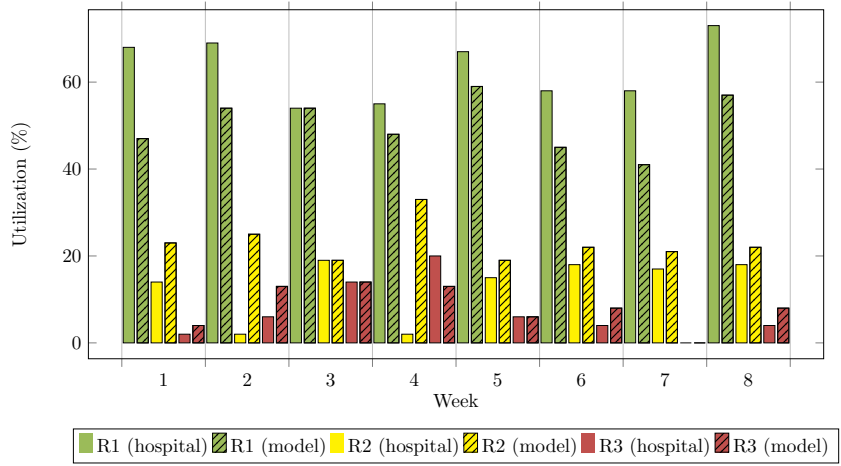
\begin{figure}[hbtp]
\caption{Utilization of units in hospital's capacity allocation policy versus optimal allocation policy} \label{v1detutil}
\centering
\resizebox{0.70\textwidth}{!}{
\begin{tikzpicture}
\centering
\begin{axis}[
	x tick label style={
		/pgf/number format/1000 sep=},
	height = 8 cm,
	width = 15 cm,
	ylabel=Utilization (\%),
	xlabel=Week,
	enlargelimits=0.05,
	legend style={at={(0.5,-0.20)},
		anchor=north,legend columns=-1},
	ybar interval=0.8,
	legend image code/.code={%
                    \draw[#1, draw=none] (0cm,-0.1cm) rectangle (0.25cm,0.25cm);
                },  
]
\addplot [style = {fill=ggreen, mark=none}]
	coordinates {(1,68) (2,69)
		 (3,54) (4,55) (5,67) (6,58) (7,58) (8,73) (9,50)};
\addplot [style = {fill=ggreen, mark=none, postaction={pattern=north east lines}}]
	coordinates {(1,47) (2,54)
		 (3,54) (4,48) (5,59) (6,45) (7,41) (8,57) (9,50)};
\addplot [style = {fill=yellow, mark=none}]
	coordinates {(1,14) (2,2)
		 (3,19) (4,2) (5,15) (6,18) (7,17) (8,18) (9,50)};
\addplot [style = {fill=yellow, mark=none, postaction={pattern=north east lines}}]
	coordinates {(1,23) (2,25)
		 (3,19) (4,33) (5,19) (6,22) (7,21) (8,22) (9,50)};
\addplot [style = {fill=rred, mark=none}]
	coordinates {(1,2) (2,6)
		 (3,14) (4,20) (5,6) (6,4) (7,0) (8,4) (9,50)};
\addplot [style = {fill=rred, mark=none, postaction={pattern=north east lines}}]
	coordinates {(1,4) (2,13)
		 (3,14) (4,13) (5,6) (6,8) (7,0) (8,8) (9,50)};
\legend{R1 (hospital) ,R1 (model) ,R2 (hospital) ,R2 (model) ,R3 (hospital) ,R3 (model)}
\end{axis}
\end{tikzpicture}
}
\end{figure}

The average daily utilization of each unit provides more insights into the clinic's current capacity allocation policy.  Figure \ref{v1detutil} presents the utilization of units in model's and hospital's allocations. We observe that for the hospital's current allocation policy, other than the third week, the standard unit's utilization is higher than that of the optimal allocation. 
In contrast, the average utilization of both isolated and quarantine units is lower than the optimal. One can argue that this makes the system more resilient to surges in suspected and infected patient numbers. We can infer that the hospital reserved more buffer capacity for the infected and suspected patients to ensure a high service level to these patients. However, this cautious approach comes at the cost of increasing overlaps among patient cohorts.\\

The results here, which hypothetically assume known demand, indicate that there might be an opportunity to make better capacity allocation decisions and reduce the number of overlapping patients in the clinic. However, in reality, weekly capacity allocation decisions must be made under significant demand uncertainty. In the next subsection, we evaluate the performance of the proposed stochastic optimization model in making capacity allocation decisions by accounting for demand uncertainty.

\subsubsection{Analysis of three-unit cohorting model's solutions} \label{s:case3b}

We now analyze the solutions of the scenario-based stochastic model, which makes capacity allocation decisions for an upcoming week under demand uncertainty to implement the three-unit cohorting strategy. As described in \S \ref{s:case2}, we generate scenarios by predicting demands for weeks six, seven and eight based on the data from the clinic. 
Table \ref{tablev1stochastic} presents the capacity allocation solution achieved by the model for these weeks under two scenario sets generated by assuming two Prediction Intervals (PI), which represent different conservatism levels to account for demand variability and the number of overlaps that would result from the capacity allocation solutions based on the realized patient demands. \\

\begin{table}[hbtp] \setlength{\tabcolsep}{4pt}
\tiny
\begin{center}
\begin{small}
\begin{tabular}{cc|lll|m{0.09\textwidth}m{0.09\textwidth}m{0.09\textwidth}l|ccc}
\cmidrule[\heavyrulewidth]{3-12}
\multicolumn{2}{c}{} & \multicolumn{7}{c}{\textbf{Three-unit cohorting }} & \multicolumn{3}{|c}{\textbf{Overlaps with }} \\
\multicolumn{2}{c}{} & \multicolumn{7}{c}{\textbf{ model's solutions (Expected overlaps)}} & \multicolumn{3}{|c}{\textbf{realized demand}} \\
\toprule
Week & PI (\%)  & $R_1$ & $R_2$ & $R_3$ & $E[O_{(1,2)x3}]$ & $E[O_{(1,2)x4}]$ & $E[O_{3x4}]$ & $Z_{3u}$ & $O_{(1,2)x3}$  & $O_{(1,2)x4}$  & $O_{3x4}$ \\
\midrule
6 & 80 & 8 & 5 & 1 & 0.0 & 0.9 & 9.3 & 1,883 & 0 & 0 & 8 \\
6 & 90 & 8 & 5 & 1 & 0.8 & 1.3 & 7.6 & 2,979 & 0 & 0 & 8 \\
7 & 80 & 8 & 5 & 1 & 1.4 & 0.5 & 10.1 & 2,998 & 0 & 0 & 0 \\
7 & 90 & 10 & 3 & 1 & 3.4 & 0.4 & 6.4 & 4,561 & 0 & 0 & 0 \\
8 & 80 & 8 & 5 & 1 & 0.0 & 0.0 & 5.1 & 561 & 24 & 5 & 9 \\
8 & 90 & 8 & 5 & 1 & 1.0 & 0.0 & 4.2 & 1,497 & 24 & 5 & 9 \\
\bottomrule
\end{tabular}
\end{small}
\end{center}
\caption{Performance of the three-unit stochastic model}
\label{tablev1stochastic}
\end{table}

Table \ref{tablev1stochastic} shows that the proposed stochastic optimization model tends to allocate a larger number of dialysis machines to the standard unit ($R_1$) that treats uninfected (Type 1 and 2) patients compared to the hospital's current allocation policy. Moreover, the capacity configuration achieved by the stochastic model leads to a lower number of overlaps for these three weeks compared to those incurred under the hospital's allocation (Table \ref{tablev1}). Specifically, in week 6, the hospital's allocation causes an overlap of uninfected and infected patients, which could be avoided with the model's solution. For week 8, the model solution can decrease the overlaps of infected and suspected patients with uninfected patients. This comes at a cost of some increase in the overlaps of suspected and infected patients, however since these overlaps are considered less critical, the model solution is an improvement over the hospital's allocation. We also observe that while week 6 and 8 capacity allocation solutions are robust to the different PIs used to develop scenarios, week 7 solutions are affected by the different set of scenarios; specifically, in week 7, the model allocates more machines to treat uninfected patients when more conservative demand realizations are accounted for. Although the number of overlaps based on realized demands in week 7 is not affected by the different capacity allocation solution, the results indicate that the prediction of demands and development of scenarios can be crucial factors affecting the performance of the capacity allocation decisions in the clinic.  

These results show that developing forecasts on the number of patients by using historical data and solving the proposed scenario-based stochastic model can help the dialysis clinic reduce the overlapping sessions of uninfected and infected patients thereby making better cohorting during the pandemic.


\subsubsection{Evaluation of different cohorting strategies} \label{s:case3c}

We next investigate the effects of cohorting strategies on the performance of the system. We first compare the solutions obtained by two- and three-unit cohorting models by considering the realized demands of the hospital for eight weeks, which are presented in Table \ref{tablev2det}. The last column of Table \ref{tablev2det} presents the difference in the objective functions of two- and three-unit models. 

\begin{table}[hbtp] \setlength{\tabcolsep}{4pt}
\begin{center}
\begin{footnotesize}
\begin{tabular}{c|ccc|cccc|cc|ccc|r}
\cmidrule[\heavyrulewidth]{2-13}
\multicolumn{1}{l}{} & \multicolumn{3}{l|}{\textbf{Optimal}}&\multicolumn{4}{l|}{\textbf{Performance with optimal}} & \multicolumn{2}{l}{\textbf{ Optimal }} & \multicolumn{3}{|l}{\textbf{Performance with optimal }} & \\
\multicolumn{1}{l}{} & \multicolumn{3}{l|}{\textbf{allocation }} &\multicolumn{4}{l|}{\textbf{allocation (known demand)}} & \multicolumn{2}{l}{\textbf{ allocation}} &\multicolumn{3}{|l}{\textbf{ allocation (known demand)}}& \\
\toprule
Week &  $R_1$ & $R_2$ & $R_3$ & $O_{(1,2)x3}$ & $O_{(1,2)x4}$ & $O_{3x4}$ & $Z_{3u}$ &$R_1$ & $R_2$ & $O_{(1,2)x3}$ & $O_{(1,2)x4}$ & $Z_{2u}$ & $\frac{Z_{3u}-Z_{2u}}{Z_{3u}}$\\
\midrule
1 &  10 & 3 & 1 & 0 & 0 & 4 & 444 &     11 & 3 & 0 & 4 & 4,042     & -810\% \\
2 &  9 & 4 & 1 & 0 & 0 & 3 & 350 &      10 & 4 & 0 & 2 & 2,050     & -486\% \\
3 &  7 & 5 & 2 & 0 & 8 & 17 & 9,752 &   9 & 5 & 0 & 11 & 11,050    & -13\% \\
4 &  8 & 3 & 3 & 0 & 0 & 17 & 1,758 &   9 & 5 & 0 & 9 & 9,048      & -415\% \\
5 &  8 & 4 & 2 & 8 & 6 & 10 & 15,050 &  10 & 4 & 0 & 13 & 13,046   & 13\% \\
6 &  9 & 4 & 1 & 0 & 0 & 5 & 546 &      8 & 6 & 0 & 2 & 2,044      & -274\% \\
7 &  10 & 4 & 0 & 0 & 0 & 0 & 36 &      10 & 4 & 0 & 0 & 36        & 0\% \\
8 &  9 & 4 & 1 & 15 & 2 & 19 & 17,950 & 10 & 4 & 8 & 13 & 21,048   & -13\% \\
\midrule
Total & - & - & - & 23 & 16 & 75 & 45,886 & - & - & 8 & 54 & 62,364  & - \\
\bottomrule
\end{tabular}
\end{footnotesize}
\end{center}
\caption{Solutions of two-unit and three unit models under deterministic demand}
\label{tablev2det}
\end{table}

As observed in Table \ref{tablev2det}, the solutions of the three-unit model dominate those of the two-unit model in seven weeks out of eight. Since infected (Type 3) and suspected (Type 4) patients can be treated in the same session in the three-unit model, the model assigns all of these patients to the last (i.e., the fourth) session of a day if possible. If the capacity of a single session is not sufficient, some patients are assigned to the third session, during which 
uninfected (Type 1 and 2) patients are usually treated as well. Therefore, by allowing overlaps in Type 3 and Type 4 patients' dialysis treatments, the three-unit model can avoid the more serious overlaps with uninfected (Type 1 and 2) patients, thus incurring smaller penalties.
  
In the two-unit case, suspected (Type 4) and infected (Type 3) patients receive dialysis in the same unit, but  suspected patients' sessions precede those of the infected ones. Consequently, since the whole unit is allocated to a single patient group when the demand from one group of patients is low, the unit capacity cannot be utilized fully.  Moreover, assigning suspected patients to sessions earlier in the day may cause increased overlaps with uninfected patients and result in worse objective values, as observed from Table \ref{tablev2det}. 


\begin{figure}[!h]
\caption{Illustration of treatment schedules for days 1 and 3 of week 5} \medskip
\label{daysweek5}
\begin{subfigure}[c]{0.48\linewidth} 
\centering 
\includegraphics[width=7.5cm]{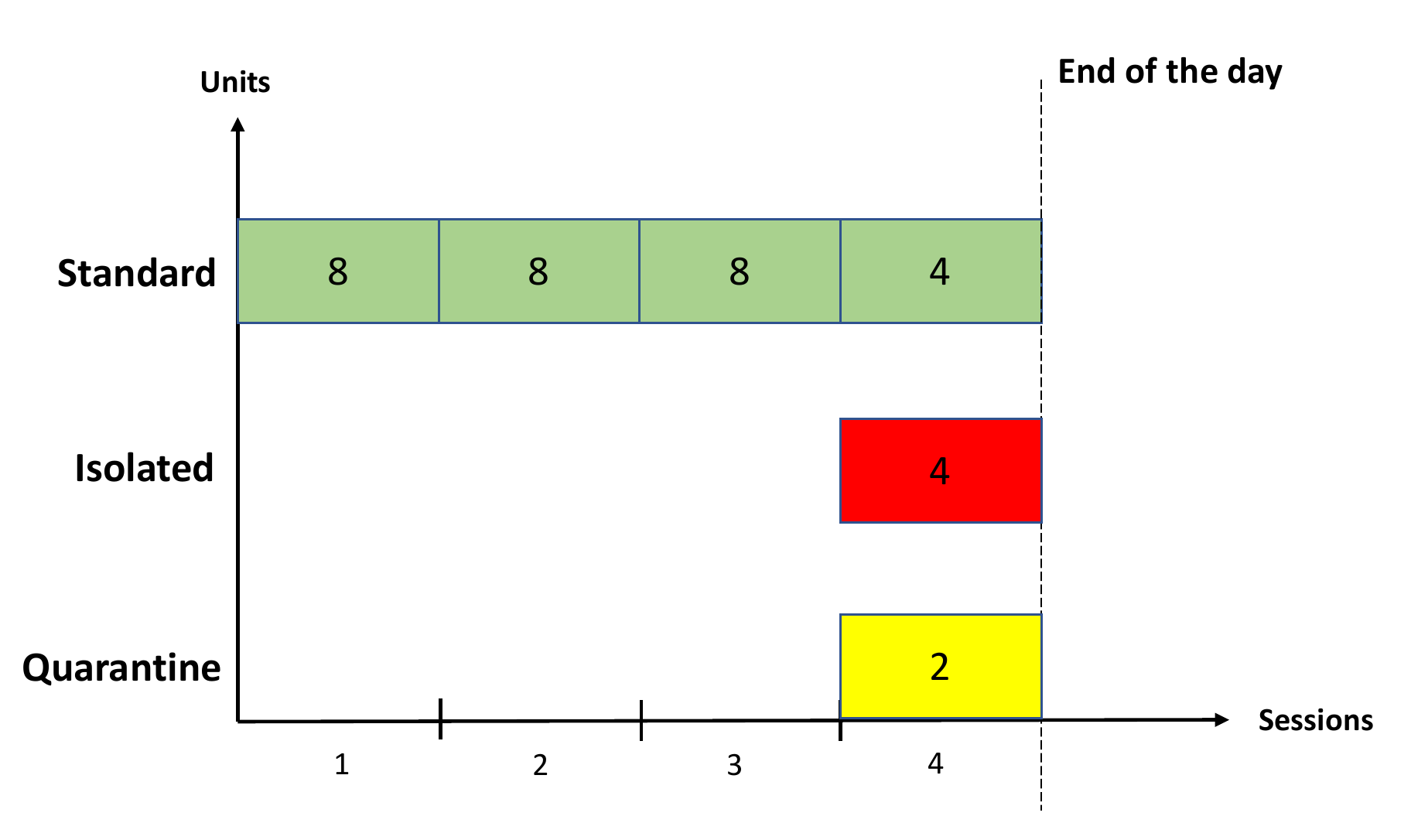} 
\caption{Day 1, three-unit model’s schedule \\ (Daily penalty = 14,612)}
\label{day13unit}
\end{subfigure}
\hfill
\begin{subfigure}[c]{0.48\linewidth} 
\centering
\includegraphics[width=7.5cm]{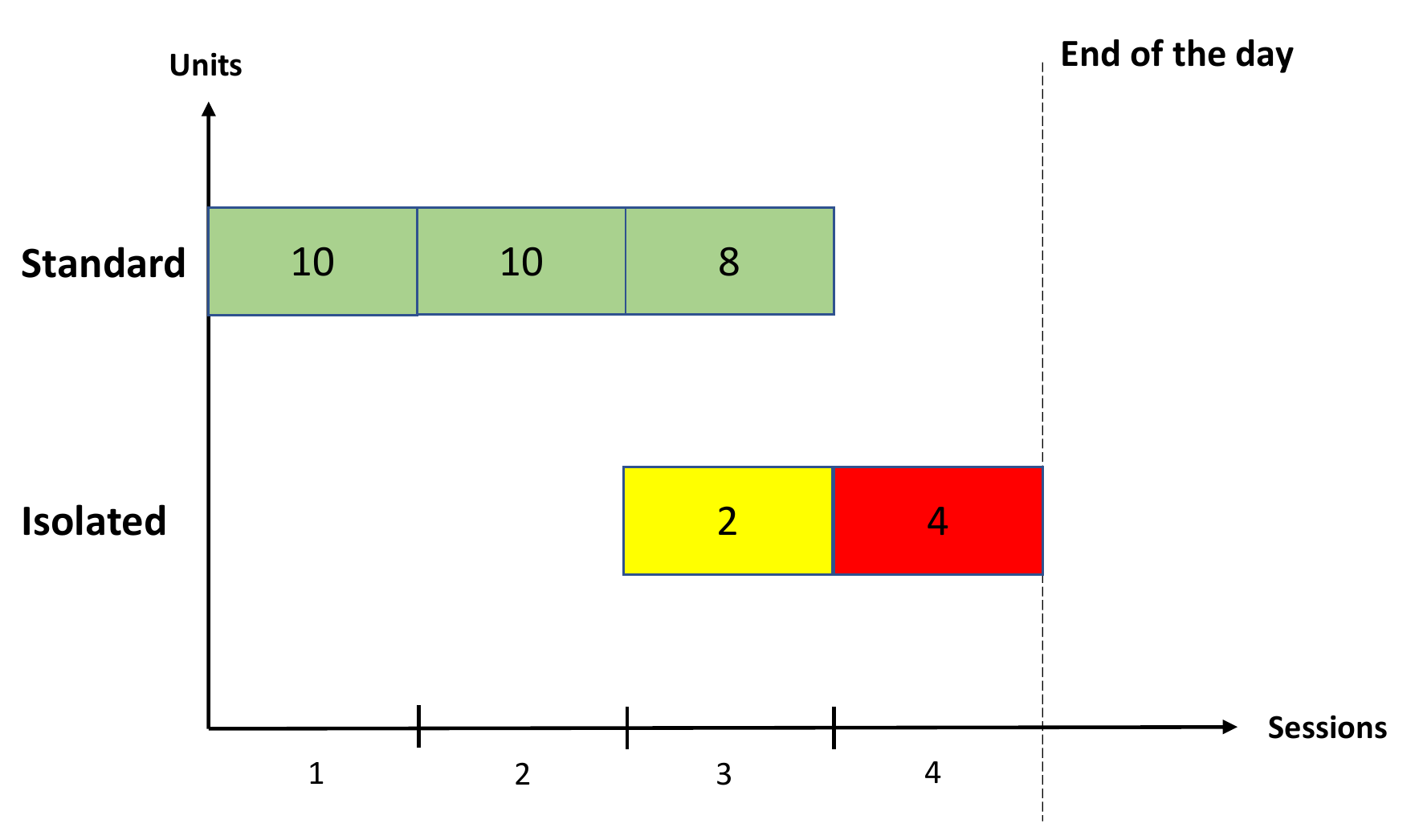} 
\caption{Day 1, two-unit model's schedule \\ (Daily penalty = 10,010)}
\label{day12unit}
\end{subfigure}

\begin{subfigure}[c]{0.48\linewidth}
\centering
\includegraphics[width=7.5cm]{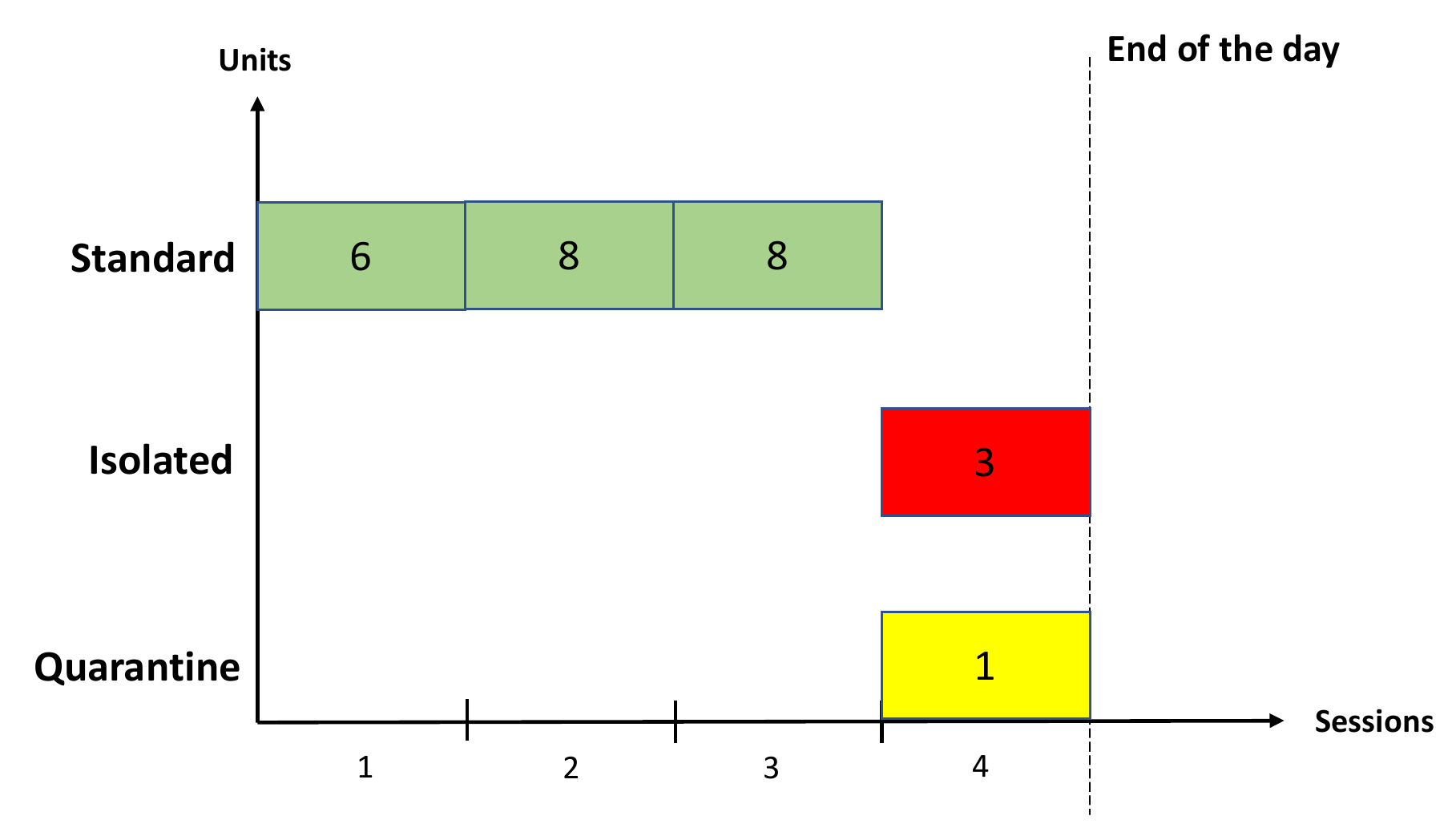} 
\caption{Day 3, three-unit model's schedule \\ (Daily penalty = 410)}
\label{day33unit}
\end{subfigure}
\hfill
\begin{subfigure}[c]{0.48\linewidth} 
\centering
\includegraphics[width=7.5cm]{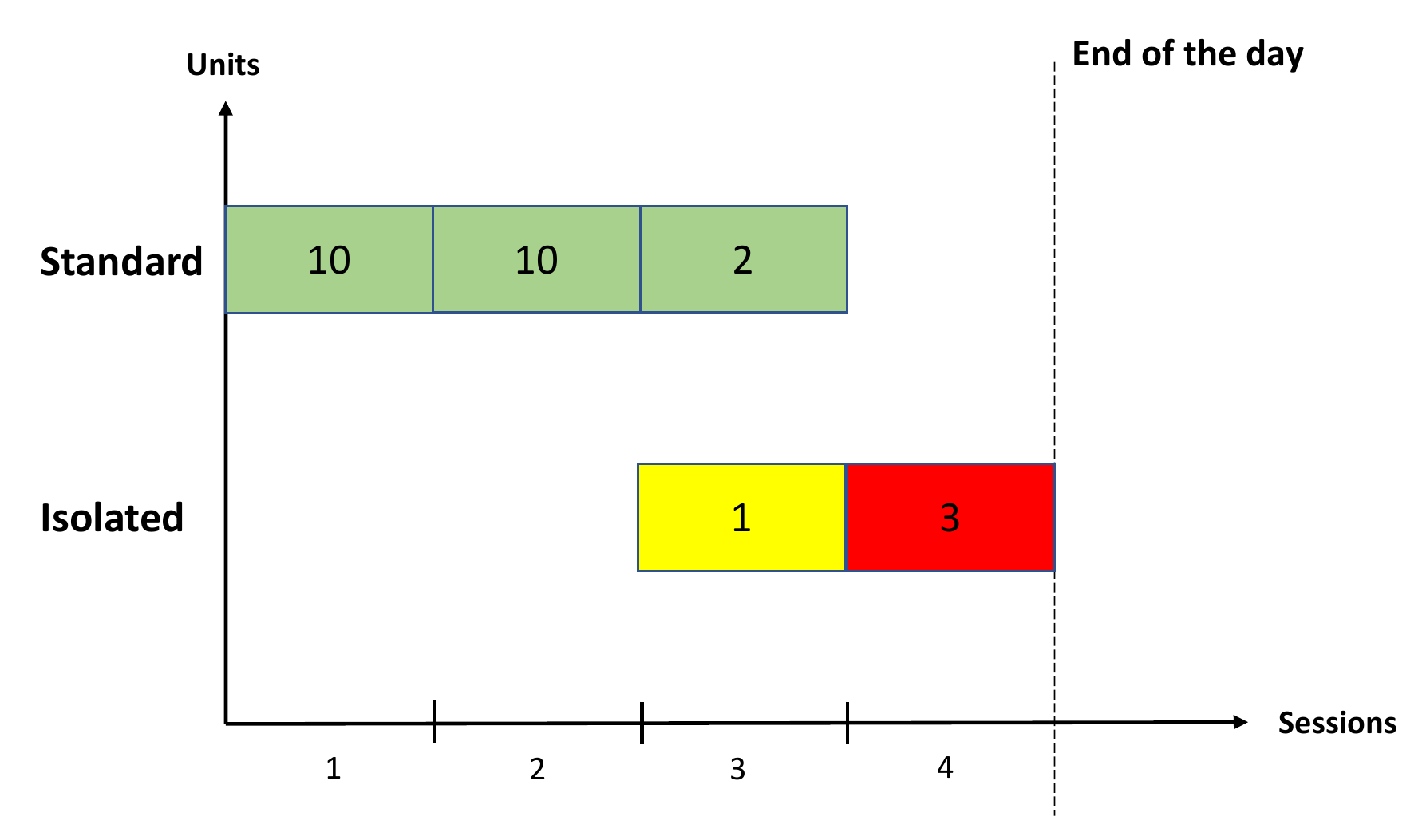}
\caption{Day 3, two-unit model's schedule \\ (Daily penalty = 3,010)}
\label{day32unit}
\end{subfigure}

\begin{subfigure}[c]{\linewidth}
\centering
\includegraphics[width=7.5cm]{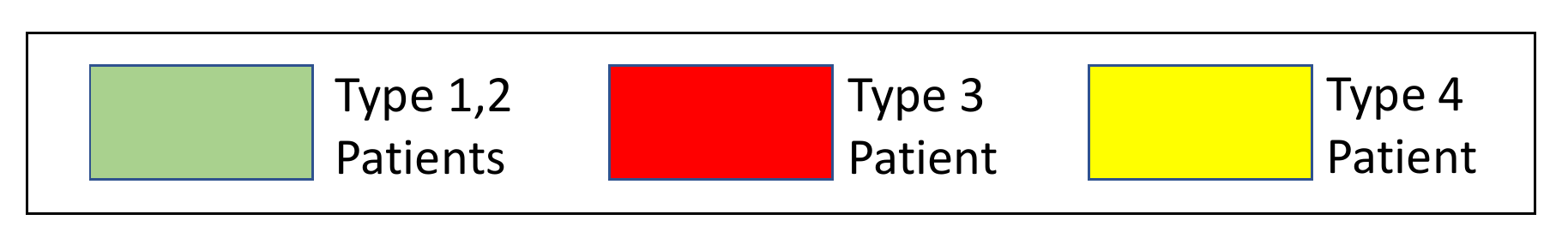} 
\end{subfigure}
\end{figure}


According to Table \ref{tablev2det}, the two-unit model slightly outperforms the three-unit model only in week 5, which is an exceptional case for the considered time horizon. Hence we examine the solutions achieved for this week in further detail. We observe that the two-unit model allocates a greater amount of resources to the standard unit compared to the three-unit model. Specifically, the optimal allocations are $R_1=10$, $R_2=4$ for the two-unit model, and $R_1=8$, $R_2=4$, $R_3=2$ for the three-unit model. In Figure \ref{daysweek5}, we present the resulting treatment schedules for two days from week 5, which include some overlaps. In week 5, the number of uninfected, infected and suspicious patient demands are 28, 4 and 2 on Monday (day 1), and 22, 1 and 3 on Wednesday (day 3). We show the number of patients assigned to each session by the alternate cohorting models in Figure \ref{daysweek5}.

As shown in Figure \ref{daysweek5}(b), on day 1, since a larger number of machines are allocated to the standard unit in the two-unit model, uninfected patients' dialysis treatments could be completed on the third session, and they only overlap with two Type 3 patients on that session. On the other hand, as shown in Figure \ref{daysweek5}(a), due to the smaller number of dialysis machines allocated to the standard unit in the three-unit case, Type 1 and 2 patients' treatments continue until the last session of the day, and a larger penalty is incurred due to several overlaps. In contrast, on day 3, the allocation made by the three-unit model results in a lower objective value (Figure \ref{daysweek5} (c)). Since Type 3 and Type 4 patients are treated in parallel in the last session, they do not overlap with any uninfected patients in this case. However, as shown in Figure \ref{daysweek5}(d), overlapping treatments had to be scheduled in the two-unit solution.


\begin{table}[h!] \setlength{\tabcolsep}{4pt}
\scriptsize
\begin{center}
\begin{small}
\begin{tabular}{cc|llm{0.09\textwidth}m{0.09\textwidth}l|cc}
\cmidrule[\heavyrulewidth]{3-9}
\multicolumn{2}{c}{} & \multicolumn{5}{c}{\textbf{Three-unit cohorting}} & \multicolumn{2}{c}{\textbf{Overlaps with}} \\
\multicolumn{2}{c}{} & \multicolumn{5}{c}{\textbf{model’s solutions (Expected overlaps)}} & \multicolumn{2}{c}{\textbf{realized demand}} \\
\toprule
Week & PI (\%) & $R_1$ & $R_2$ & $E[O_{(1,2)x3}]$ & $E[O_{(1,2)x4}]$ & $Z_{2u}$ & $O_{(1,2)x3}$ & $O_{(1,2)x4}$ \\ 
\midrule
6 & 80 & 9 & 5 & 0.0 & 4.7 & 4713 & 0 & 7 \\ 
6 & 90 & 9 & 5 & 0.0 & 4.1 & 4112 & 0 & 7 \\ 
7 & 80 & 9 & 5 & 0.9 & 5.9 & 6815 & 0 & 0 \\ 
7 & 90 & 8 & 6 & 0.0 & 7.2 & 7214 & 0 & 0 \\ 
8 & 80 & 9 & 5 & 0.0 & 2.2 & 2211 & 16 & 14 \\
8 & 90 & 9 & 5 & 0.7 & 3.1 & 3811 & 16 & 14 \\
\bottomrule
\end{tabular}
\end{small}
\end{center}
\caption{Performance of the two-unit stochastic model}
\label{tablev2stochastic}
\end{table}


We also evaluate the performance of the two-unit cohorting strategy by solving the two-unit scenario-based stochastic model, in which the uncertainties inpatient demands are accounted for while deciding on the number of machines for an upcoming week. Table \ref{tablev2stochastic} presents the results of the two-unit model, which are compared with those obtained by the three-unit model for the same instances (Table \ref{tablev1stochastic}). 
We observe that the three-unit cohorting strategy results in lower expected penalties than those of the two-unit cohorting strategy.
    
These results show that determining the best cohorting policy and the number of dialysis machines to allocate to different units can be challenging, even under the availability of perfect demand information. The demands for each patient type, which can be highly uncertain, can significantly affect the performance of cohorting strategies; and it is hard to foresee which cohorting strategy would perform better without solving the stochastic optimization models proposed in this paper. 
The complexity of the cohorting, capacity allocation and treatment scheduling decisions, which increase under demand uncertainty, imply the need for analytical approaches to effectively manage dialysis treatments of various patient types during a pandemic. Using the methods proposed in this paper, dialysis clinics can evaluate the performance of alternative cohorting policies and obtain the best allocation decision based on the forecasts for the upcoming time horizon.

\section{Conclusion} \label{s:conclusion}

The COVID-19 pandemic has shown that effective planning of health care resources is crucial. Moreover, treating vulnerable groups such as people with chronic diseases needs special attention to prevent infection spread in health facilities. In this paper, we have focused on the operations of a hemodialysis clinic in a major pandemic hospital in Turkey, which is providing dialysis treatment to multiple patient groups. While the hospital has taken effective precautions to manage the treatment of different types of patients through cohorting, they lack the analytical tools to make cohorting plans and to evaluate the effectiveness of different strategies. We present a two-stage stochastic programming modeling approach to help the clinic make more effective capacity planning and treatment scheduling decisions. We show with real data collected during the COVID-19 pandemic that the clinic can make use of such an operations research based tool to mitigate the infection transmission risk at the hospital by decreasing overlapping sessions among infected and uninfected patient groups. Moreover, we present results that compare two alternative cohorting strategies based on the hospital's data. We show that while three-unit cohorting is generally more effective than two-unit cohorting to reduce overlaps among patient cohorts in our case study, the performance of each strategy ican be highly dependent on the relative number of each patient type. The proposed analytical methods can support the hospital to cohort patients during a pandemic effectively and manage scarce health resources efficiently. 

Although the hospital was able to treat all patients with the available dialysis resources since the beginning of the pandemic, albeit incurring some overlapping sessions, the proposed models can also be used to assess when it might be impossible to treat all patients in the clinic with the existing capacity should the demands from COVID-19 patients increase due to another pandemic wave. In such a case, the model results can be presented to policymakers that can consider making arrangements to transfer chronic dialysis patients to other private facilities in emergency cases, which is not currently allowed. Moreover, the non-medical operational solution proposed in this paper can be used to support cohorting, capacity planning and treatment scheduling decisions in different hospital units besides hemodialysis such as to manage the treatments of chemotherapy, radiotherapy and physiotherapy patients.

Given the scarce research attention on the planning and coordination of health care resources to treat chronic patients during extreme events such as natural disasters and pandemics, several future research directions exist. For example, the proposed models can be extended to develop tools for managing other types of chronic patients (such as chemotherapy patients) who may present different aspects and constraints. While dialysis treatment session durations are fixed, chemotherapy or physical therapy treatment durations may be different for each patient. This makes the capacity per day dependent on the patients scheduled, and the scheduling problem becomes more challenging due to the additional complexity of matching different treatment types. Another future research direction can be to incorporate uncertainty in patient arrivals during the day. We assumed that the demands from all patient types are known at the beginning of the day while scheduling treatments to minimize overlaps. Arrivals of patients with urgent treatment needs can be accommodated by reserving capacity for these patients, which is a common practice in appointment systems. Our models can be extended in future research to account for such dynamics. Finally, different cohorting strategies can be analyzed with queueing models to obtain insights on their effectiveness in different environments.



\bibliography{references}

\section*{Acknowledgements}
This research has been supported by AXA Award Grant from AXA Research Fund.
\section{Appendices}
\subsection*{A. Case Data} \label{Appendix1}
\begin{table}[hp!]
  \centering
    \begin{adjustbox}{width=\textwidth}
\begin{tabular}{|l|llll|llll|llll|llll|llll|llll|llll|}
\hline
\textbf{}            & \multicolumn{4}{c}{\textbf{Monday}} & \multicolumn{4}{c}{\textbf{Tuesday}} & \multicolumn{4}{c}{\textbf{Wednesday}} & \multicolumn{4}{c}{\textbf{Thursday}} & \multicolumn{4}{c}{\textbf{Friday}} & \multicolumn{4}{c}{\textbf{Saturday}} & \multicolumn{4}{c|}{\textbf{Sunday}} \\
\hline
\textbf{Patient type:}               & \textbf{1}       & \textbf{2}        & \textbf{3}      & \textbf{4}      & \textbf{1}       & \textbf{2}        & \textbf{3}      & \textbf{4} & \textbf{1}       & \textbf{2}        & \textbf{3}      & \textbf{4} & \textbf{1}       & \textbf{2}        & \textbf{3}      & \textbf{4} & \textbf{1}       & \textbf{2}        & \textbf{3}      & \textbf{4} & \textbf{1}       & \textbf{2}        & \textbf{3}      & \textbf{4} & \textbf{1}       & \textbf{2}        & \textbf{3}      & \textbf{4}     \\
\hline
\textbf{Week 1}             & 6       & 12      & 2      & 0      & 11      & 8       & 2       & 0      & 8       & 12       & 3       & 0       & 13      & 8       & 2       & 0       & 13      & 12      & 3      & 1      & 3       & 8       & 5       & 0       & 1       & 0       & 0      & 0      \\
\textbf{Week 2}             & 13      & 12      & 4      & 0      & 8       & 8       & 4       & 1      & 11      & 12       & 4       & 0       & 7       & 8       & 6       & 0       & 9       & 12      & 2      & 1      & 8       & 8       & 4       & 1       & 1       & 0       & 2      & 0      \\
\textbf{Week 3}             & 7       & 12      & 2      & 1      & 2       & 8       & 5       & 1      & 7       & 12       & 5       & 2       & 4       & 8       & 3       & 0       & 9       & 12      & 5      & 3      & 2       & 8       & 3       & 0       & 4       & 0       & 2      & 0      \\
\textbf{Week 4}             & 12      & 12      & 2      & 3      & 6       & 8       & 5       & 1      & 3       & 12       & 4       & 3       & 4       & 8       & 5       & 0       & 4       & 12      & 4      & 2      & 4       & 8       & 4       & 1       & 3       & 0       & 1      & 0      \\
\textbf{Week 5}             & 16      & 12      & 4      & 2      & 7       & 8       & 2       & 0      & 10      & 12       & 3       & 1       & 7       & 8       & 2       & 0       & 6       & 12      & 3      & 0      & 8       & 8       & 5       & 0       & 4       & 0       & 3      & 0      \\
\textbf{Week 6}             & 8       & 12      & 1      & 0      & 9       & 8       & 5       & 1      & 3       & 12       & 3       & 0       & 7       & 8       & 3       & 0       & 4       & 12      & 4      & 0      & 7       & 8       & 6       & 1       & 2       & 0       & 2      & 0      \\
\textbf{Week 7}             & 6       & 12      & 4      & 0      & 8       & 8       & 4       & 0      & 3       & 12       & 3       & 0       & 10      & 8       & 4       & 0       & 5       & 12      & 3      & 0      & 7       & 8       & 3       & 0       & 3       & 0       & 0      & 0      \\
\textbf{Week 8}             & 13      & 12      & 7      & 0      & 10      & 8       & 3       & 0      & 11      & 12       & 3       & 1       & 5       & 8       & 3       & 0       & 16      & 12      & 4      & 1      & 9       & 8       & 2       & 0       & 6       & 0       & 0      & 0      \\
\hline
\textbf{Average}            & 10.1   & 12.0   & 3.3   & 0.8   & 7.6    & 8.0   & 3.8   & 0.5   & 7.0    & 12.0    & 3.5    & 0.9    & 7.1    & 8.0    & 3.5    & 0.0    & 8.3    & 12.0   & 3.5   & 1.0   & 6.0    & 8.0  & 4.0    & 0.4    & 3.0    & 0.0    & 1.3   & 0.0   \\
\textbf{Minimum}            & 6       & 12      & 1      & 0      & 2       & 8       & 2       & 0      & 3       & 12       & 3       & 0       & 4       & 8       & 2       & 0       & 4       & 12      & 2      & 0      & 2       & 8       & 2       & 0       & 1       & 0       & 0      & 0      \\
\textbf{Maximum}            & 16      & 12      & 7      & 3      & 11      & 8       & 5       & 1      & 11      & 12       & 5       & 3       & 13      & 8       & 6       & 0       & 16      & 12      & 5      & 3      & 9       & 8       & 6       & 1       & 6       & 0       & 3      & 0      \\
\textbf{Std. Dev.} & 3.8    & 0.0    & 1.9   & 1.2   & 2.8    & 0.0    & 1.3    & 0.5   & 3.6    & 0.0     & 0.8    & 1.1    & 3.1    & 0.0    & 1.4    & 0.0    & 4.40    & 0.00    & 0.9   & 1.1   & 2.6    & 0.0    & 1.3    & 0.5    & 1.7    & 0.0    & 1.2   & 0.0  
 \\ 
\hline
\end{tabular}
\end{adjustbox}
\caption{Raw data and descriptive statistics on daily patients demands (November 2-December 28)} \medskip
  \label{rawdata}
\end{table}


\subsection*{B. Discrete Distributions for Patient Types} \label{Appendix2}

\label{Appendix2}
\begin{table}[H] \setlength{\tabcolsep}{4pt}
\begin{center}
\begin{footnotesize}
\medskip
\begin{tabular}{cc|lllllll}
\toprule
\multicolumn{2}{c}{\textbf{}} & \multicolumn{7}{c}{\textbf{P(X = x)}} \\
\toprule
\textbf{Patient Type} & \textbf{PI (\%)} & \textbf{0} & \textbf{1} & \textbf{2} & \textbf{3} & \textbf{4} & \textbf{5} & \textbf{6} \\
\midrule
Type 3 & 80 & & & 0.189 & 0.302 & 0.302 & 0.207 & \\
Type 3 & 90 & & 0.021 & 0.236 & 0.236 & 0.236 & 0.236 & 0.035 \\
Type 4 & 80 & 0.688 & 0.312 & & & & & \\
Type 4 & 90 & 0.647 & 0.353 & & & & & \\
\bottomrule
\end{tabular}
\end{footnotesize}
\end{center}
\caption{Discrete probability distributions obtained for Type 3 and 4 patients for week 8}
\label{table:discdist2}
\end{table}

\end{document}